\documentclass[11pt]{article}
\usepackage[margin=1in]{geometry}

\usepackage{amsmath,amssymb,amsthm}
\usepackage{proof} 
\usepackage{hyperref}
\usepackage[T1]{fontenc}
\usepackage{enumerate}
\usepackage{mathtools}
\usepackage{xcolor}
\usepackage{array}
\usepackage{booktabs}

\hypersetup{
	colorlinks=true,
	linkcolor=blue!60!black,
	citecolor=green!50!black,
	urlcolor=blue!60!black
}

\newtheorem{theorem}{Theorem}[section]
\newtheorem{lemma}[theorem]{Lemma}
\newtheorem{corollary}[theorem]{Corollary}

\newtheorem{definition}[theorem]{Definition}
\newtheorem{remark}[theorem]{Remark}

\newcommand{\QRC}{\mathrm{QRC}_1}
\newcommand{\GQRC}{\mathrm{G3QRC}_1}

\newcommand{\GL}{\mathrm{GL}}

\newcommand{\PA}{\mathrm{PA}}
\newcommand{\HA}{\mathrm{HA}}
\newcommand{\seq}{\Rightarrow}
\newcommand{\cpx}[1]{|{#1}|}
\newcommand{\rname}[1]{\mathsf{#1}}
\newcommand{\rId}{\rname{Id}}
\newcommand{\rTopR}{\rname{\top R}}
\newcommand{\rAndL}{\rname{\wedge L}}
\newcommand{\rAndR}{\rname{\wedge R}}
\newcommand{\rAllL}{\rname{\forall L}}
\newcommand{\rAllR}{\rname{\forall R}}
\newcommand{\rDiamK}{\rname{\Diamond K}}
\newcommand{\rDiamFour}{\rname{\Diamond 4}}
\newcommand{\rCut}{\rname{Cut}}
\newcommand{\fv}{\mathrm{fv}}
\newcommand{\vocab}{\mathrm{vocab}}
\newcommand{\consts}{\mathrm{const}}
\newcommand{\Hderiv}{\vdash}

\begin{document}
	
	\title{A Sequent Calculus for the Quantified Reflection
		Calculus $\QRC$: Cut-Elimination, Interpolation,
		and NP-Completeness}
	
	\author{Mitar Nedeljkovi\'c}
	\date{}
	\maketitle
	
	\begin{abstract}
			We introduce $\GQRC$, the first Gentzen-style sequent calculus for Borges and Joosten's Quantified Reflection Calculus ($\QRC$). As a strictly positive, first-order modal logic, $\QRC$ crucially escapes Vardanyan's $\Pi^0_2$-completeness theorem for quantified provability logic, while remaining arithmetically complete with respect to Peano and Heyting Arithmetic. The central innovation of $\GQRC$ is a pair of \emph{context-erasing} modal rules that simultaneously account for the monotonicity and transitivity axioms by deliberately discarding parametric context across modal boundaries. This structural bottleneck allows for a clean proof of cut-admissibility via a standard lexicographic induction, circumventing the complex loop-checking typically required for provability logics. Cut admissibility yields an instantiated subformula property and supports a direct canonical proof of the finite model property. We streamline the adequate Kripke semantics by interpreting constants globally. We also establish a split-context Craig interpolation theorem and a contextual form of left uniform interpolation, uniformly eliminating private predicate symbols, constants, and free parameters. Finally, by leveraging the rigid structural constraints of the context-erasing rules and the absence of left-branching connectives, we prove that, over the standard fixed countable signatures specified in Section~\ref{sec:complexity}, the decision problem for $\QRC$ is $\mathbf{NP}$-complete, situating this quantified modal logic in a comparatively low nondeterministic complexity class.
		\end{abstract}
	\medskip
	
	
	\bigskip
	
	\section{Introduction}
	\label{sec:intro}
	
	The study of provability logic has historically been characterized by a stark dichotomy between the propositional and the quantified realms. At the propositional level, G\"{o}del-L\"{o}b logic ($\GL$) provides a pristine, decidable axiomatization of provability in Peano Arithmetic ($\PA$), famously validated by Solovay's completeness theorem \cite{Solovay1976}. However, extending this elegant framework to the first-order level precipitates a collapse into intractability. As demonstrated by Vardanyan \cite{Vardanyan1986}, the quantified provability logic of $\PA$ ($\mathrm{QPL}(\PA)$) is $\Pi^0_2$-complete, rendering it not even recursively axiomatizable. This barrier is remarkably robust, extending to a wide class of arithmetical theories, a phenomenon Visser and de Jonge aptly termed ``No escape from Vardanyan's theorem'' \cite{VisserDeJonge2006}.
	
	A profound escape route was eventually discovered through the paradigm of \emph{strictly positive} logics. Developed by Beklemishev (\cite{Beklemishev2004}, \cite{Beklemishev2014}) for the study of uniform reflection principles and ordinal analysis via the Reflection Calculus ($\mathrm{RC}$), strictly positive logics deliberately omit negation and implication. The Vardanyan-style hardness construction combines arithmetical coding of computation with diagonal, self-referential machinery; the strictly positive language blocks this route while preserving the expressive power needed to formalize relative consistency. Building on this foundation, de Almeida Borges and Joosten \cite{BorgesJoosten2020} successfully lifted the strictly positive paradigm to the quantified level, introducing the Quantified Reflection Calculus ($\QRC$). By interpreting predicate symbols as parameterized families of theories, they bypassed Vardanyan's theorem entirely. $\QRC$ has since proven to be a remarkably fertile system, yielding arithmetical completeness with respect to both $\PA$ \cite{BorgesJoosten2023} and Heyting Arithmetic ($\HA$) \cite{BorgesJoosten2024}, and has been formalized in part within the Coq proof assistant \cite{Borges2022Coq}.
	
	Despite these immense semantic and arithmetical triumphs, $\QRC$ has thus far lacked a native Gentzen-style proof theory. The original formulation of $\QRC$ is a Hilbert-style calculus. While perfectly suited for algebraic and arithmetical proofs, Hilbert systems inherently provide little structural control over derivations. Without an analytic, connective-by-connective decomposition of formulas, one cannot easily employ root-first proof search, establish subformula properties, or extract precise computational complexity bounds. 
	
	The present paper fills this structural gap by introducing $\GQRC$, the first analytic sequent calculus for $\QRC$. Constructing a well-behaved Gentzen calculus for provability logics is notoriously difficult due to the transitivity axiom. In standard sequent calculi, transitivity rules duplicate principal formulas into the premises, paralyzing the inductive measure required to prove cut-elimination, a challenge famously navigated by Valentini for $\GL$ \cite{Valentini1983} and extensively studied in structural proof theory \cite{Poggiolesi2010}. 
	
	We resolve this obstacle by introducing a pair of \emph{context-erasing} modal rules, $\rDiamK$ and $\rDiamFour$. These rules force all side-formulas to be permanently discarded when traversing a modal boundary. This mechanism makes the modal commutative case immediate and forces the lexicographic cut measure to decrease: usually the cut-formula complexity decreases, and in the remaining cases the sum of the premise heights decreases. This allows for a remarkably clean, standard lexicographic induction for cut-admissibility.
	
	The analytic power of $\GQRC$ unlocks a comprehensive suite of new results for $\QRC$, vastly expanding upon the initial semantic findings. By leveraging the cut-elimination theorem, we establish the subformula property. Strengthening the single-antecedent interpolation theorem of de Almeida Borges and Joosten \cite{BorgesJoosten2024}, we prove a split-context Craig interpolation theorem and contextual left uniform interpolation. The interpolation construction is a semantic projection that uniformly forgets private predicate symbols, individual constants, and free parameters. Furthermore, we streamline the original adequate Kripke models \cite{BorgesJoosten2020} by evaluating constants globally, providing a direct, self-contained proof of the finite model property via root-first saturation. Finally, the structural rigidity of $\GQRC$ allows us to precisely pinpoint the computational complexity of the logic, solving a crucial open dimension of $\QRC$. While strictly positive fragments typically lower complexity, the interaction between first-order quantifiers and modalities often resurrects $\mathbf{EXPTIME}$ or undecidable behaviors. By analyzing the derivation space of $\GQRC$, we demonstrate that the context-erasing rules bound the number of modal inferences along each proof branch and prevent unselected antecedent formulas from propagating across modal inferences. We prove that, over the standard fixed countable signatures specified in Section~\ref{sec:complexity}, the decision problem for $\QRC$ is $\mathbf{NP}$-complete, situating this quantified modal logic in a comparatively low nondeterministic complexity class.
	
	After the submission of the present manuscript on 2 June 2026, Perini Brogi and Joosten posted an independently developed preprint on 31 October 2026 \cite{PeriniBrogiJoosten2026}, establishing cut admissibility and uniform interpolation for $\QRC$ using a single-antecedent sequent calculus, in contrast to the multiset-based calculus developed here.
	
	The paper is structured as follows. Section~\ref{sec:qrc} recalls the syntax, Hilbert calculus, and semantics of $\QRC$. Section~\ref{sec:calculus} introduces the $\GQRC$ calculus. Section~\ref{sec:structural} establishes its foundational structural properties, culminating in the admissibility of contraction. Section~\ref{sec:cut} presents the proof of cut-admissibility. Section~\ref{sec:completeness} proves soundness and completeness with respect to both the Hilbert system and the semantics, yielding the finite model property. Section~\ref{sec:consequences} derives the subformula property and the interpolation theorems. Finally, Section~\ref{sec:complexity} establishes the $\mathbf{NP}$-completeness bounds.
	
	\section{The Logic \texorpdfstring{$\QRC$}{QRC1}}
	\label{sec:qrc}
	
	\subsection{Language and Hilbert axiomatization\\}
	\label{subsec:hilbert}
	
	We write $\mathsf{Var}=\{x_0,x_1,\ldots\}$ for a fixed countable set of \emph{variables}, and fix an at most countable \emph{signature} $\Sigma$ consisting of at most countably many constants and relation symbols with specified finite arities, without any function symbols. A \emph{term} is either a variable or a constant.
	
	\begin{definition}[Strictly positive formulas]
		The set $\mathcal{L}^+$ of \emph{strictly positive formulas} over $\Sigma$ is generated by the grammar:
		\[
		\varphi \;::=\; \top \;\mid\; S(t_1,\ldots,t_n) \;\mid\; \varphi \wedge \psi \;\mid\; \Diamond\varphi \;\mid\; \forall x\,\varphi
		\]
		where $S$ ranges over $n$-ary relation symbols of $\Sigma$, $t_i$ over terms, and $x$ over variables.	There is no negation, disjunction, implication, or box modality. The free variables $\fv(\varphi)$ are defined in the standard way; for terms, $\fv(x)=\{x\}$ and $\fv(c)=\varnothing$. Let $\mathrm{Var}(\varphi)$ denote the set of all variable symbols occurring in $\varphi$, free or bound, and define $\mathrm{Var}(t)$ analogously for terms. The \emph{complexity} $\cpx{\varphi}$ is the total number of logical symbols ($\wedge$, $\Diamond$, $\forall$) in $\varphi$; atomic formulas and $\top$ have complexity $0$. The \emph{modal depth} $\mathrm{md}(\varphi)$ is defined recursively by
		\[
		\mathrm{md}(\top)=\mathrm{md}(P)=0,\qquad
		\mathrm{md}(A\wedge B)=\max\{\mathrm{md}(A),\mathrm{md}(B)\},
		\]
		\[
		\mathrm{md}(\forall xA)=\mathrm{md}(A),\qquad
		\mathrm{md}(\Diamond A)=\mathrm{md}(A)+1.
		\]
	\end{definition}
	
	\begin{definition}[Substitution]
		\label{subsdef}
		In the fresh-variable clauses below, the fresh variable is the least member of the fixed enumeration $x_0,x_1,\ldots$ satisfying the displayed condition. For a term $t$, a variable $x$, and a formula $\varphi$, the \emph{substitution} $\varphi[t/x]$ is defined recursively by:
		\begin{align*}
			\top[t/x] &= \top, \\
			S(t_1,\dots,t_n)[t/x] &= S(t_1[t/x],\dots,t_n[t/x]), \\
			(\varphi\wedge\psi)[t/x] &= \varphi[t/x] \wedge \psi[t/x], \\
			(\Diamond\varphi)[t/x] &= \Diamond(\varphi[t/x]), \\
			(\forall z\,\varphi)[t/x] &=
			\begin{cases}
				\forall z\,\varphi & \text{if } z = x, \\
				\forall z\,(\varphi[t/x]) & \text{if } z \neq x
				\text{ and } z \notin \fv(t), \\
				\forall z'\,(\varphi[z'/z][t/x]) & \text{otherwise,}
			\end{cases}
		\end{align*}
		where in the last clause $z'$ is the least variable such that $z'\notin\mathrm{Var}(\varphi)\cup\mathrm{Var}(t)\cup\{x\}$, and $t_1[t/x],\dots$ denote the obvious substitution on terms. The condition ``$t$ is free for $x$ in $\varphi$'' has its usual meaning: no free occurrence of $x$ in $\varphi$ lies in the scope of a quantifier $\forall z$ with $z\in\fv(t)$.
	\end{definition}
	
	Substitution of a term for a constant, written $\varphi[t/c]$, is also capture-avoiding. On terms we set $c[t/c]=t$, $c'[t/c]=c'$ for $c'\neq c$, and $x[t/c]=x$. The clauses for $\top$, atoms, conjunction, and $\Diamond$ are as above. For a quantifier,
	\[
	(\forall z\,\varphi)[t/c]=
	\begin{cases}
		\forall z\,(\varphi[t/c])&\text{if }z\notin\mathrm{Var}(t),\\
		\forall z'\,(\varphi[z'/z][t/c])&\text{if }z\in\mathrm{Var}(t),
	\end{cases}
	\]
	where $z'$ is the least variable such that $z'\notin\mathrm{Var}(\varphi)\cup\mathrm{Var}(t)$. 
	
	\begin{definition}[$\alpha$-Equivalence]
		\label{alphadef}
		Two formulas $\varphi$ and $\psi$ are \emph{$\alpha$-equivalent},
		denoted $\varphi\equiv_\alpha\psi$, if they differ only by the
		capture-avoiding renaming of bound variables. Formally,
		$\equiv_\alpha$ is defined recursively on formula complexity by:
		\begin{enumerate}
			\item Atomic formulas are $\alpha$-equivalent iff they are
			syntactically identical, and $\top\equiv_\alpha\top$.
			\item $A\wedge B\equiv_\alpha A'\wedge B'$ iff
			$A\equiv_\alpha A'$ and $B\equiv_\alpha B'$.
			\item $\Diamond A\equiv_\alpha\Diamond A'$ iff
			$A\equiv_\alpha A'$.
			\item $\forall x\gamma\equiv_\alpha\forall y\varphi$ iff, for
			every completely fresh variable
			\[
			z\notin
			\mathrm{Var}(\forall x\gamma)\cup
			\mathrm{Var}(\forall y\varphi),
			\]
			one has
			\[
			\gamma[z/x]\equiv_\alpha\varphi[z/y].
			\]
		\end{enumerate}
		Formulas with different outermost constructors are not
		$\alpha$-equivalent.
	\end{definition}

	\begin{lemma}[$\alpha$-Stability of Substitution]
		\label{lem:alpha-subst}
		If $A\equiv_\alpha A'$, then, for every term $t$, variable $x$, and constant $c$,
		\[
		A[t/x]\equiv_\alpha A'[t/x]
		\qquad\text{and}\qquad
		A[t/c]\equiv_\alpha A'[t/c].
		\]
		Moreover, for every finite set of variables $V$, each $\alpha$-equivalence class has a representative all of whose bound variables avoid $V$.
	\end{lemma}
	\begin{proof}
		The two substitution claims are proved simultaneously by induction on the common formula structure. The atomic, conjunction, and modal cases are immediate. In the quantifier case, first rename the displayed binders in both formulas to one variable fresh from the two bodies, from $\mathrm{Var}(t)$, and from the variable being substituted; the capture-avoiding quantifier clauses then pass both substitutions through the same binder, and the induction hypothesis applies to the bodies. The final assertion follows by induction on formula complexity, choosing a fresh binder at each quantifier from the infinite variable supply.
	\end{proof}
	
	A \emph{finite simultaneous substitution} is a map
	$\sigma:\mathsf{Var}\to\{\text{terms}\}$ with finite support
	\[
	\operatorname{supp}(\sigma)=\{x:\sigma(x)\neq x\}.
	\]
	For a term $s$, put $s\sigma=\sigma(s)$ when $s$ is a variable and
	$s\sigma=s$ when $s$ is a constant. To define $A\sigma$, first choose an
	$\alpha$-equivalent representative whose bound variables avoid
	$\operatorname{supp}(\sigma)$ and all variables occurring in the
	nonidentity values of $\sigma$, and then replace all free variables
	simultaneously. A structural induction shows that different standardized
	representatives give $\alpha$-equivalent results, so $A\sigma$ is well
	defined up to $\alpha$-equivalence.
	
	\begin{lemma}[Substitution Composition]
		\label{lem:subst-comp}
		Let $\sigma$ be a finite simultaneous substitution and let $t$ be a term.
		Define $\tau$ by
		\[
		\tau(x)=t\sigma,
		\qquad
		\tau(z)=\sigma(z)\quad(z\neq x).
		\]
		After standardizing the bound variables of $\psi$ apart from
		\[
		\{x\}\cup\mathrm{Var}(t)\cup\operatorname{supp}(\sigma)
		\cup\bigcup_{z\in\operatorname{supp}(\sigma)}
		\mathrm{Var}(\sigma(z)),
		\]
		one has
		\[
		(\psi[t/x])\sigma\equiv_\alpha\psi\tau.
		\]
	\end{lemma}
	\begin{proof}
		By structural induction on $\psi$. The atomic, conjunction, and modal cases are immediate. In the quantifier case both substitutions pass through the same standardized-apart binder, and the induction hypothesis applies to the body.
	\end{proof}
	
	\begin{definition}[$\QRC$: the Quantified Reflection Calculus with one modality]
		\label{def:QRC}
		Fix a signature $\Sigma$. The Hilbert system $\QRC$ derives \emph{statements} of the form $\varphi \Hderiv \psi$ where $\varphi, \psi \in \mathcal{L}^+$. The axioms and rules are the following (letting $\varphi,\psi,\gamma$ range over $\mathcal{L}^+$) \cite[p.~15]{BorgesJoosten2020}:
		
		\smallskip
		\begin{center}
			\begin{tabular}{r@{\quad}l@{\qquad}r@{\quad}l}
				(i) & $\varphi \Hderiv \top$ and $\varphi \Hderiv \varphi$; & 
				(viii)& if $\varphi \Hderiv \psi$, then $\varphi \Hderiv \forall x \psi$ \\
				& & & ($x \notin \fv(\varphi)$); \\[1ex]
				(ii) & $\varphi \wedge \psi \Hderiv \varphi$ and $\varphi \wedge \psi \Hderiv \psi$; & 
				(ix) & if $\varphi[t/x] \Hderiv \psi$ then $\forall x \varphi \Hderiv \psi$ \\
				& & & ($t$ free for $x$ in $\varphi$); \\[1ex]
				(iii) & if $\varphi \Hderiv \psi$ and $\varphi \Hderiv \gamma$, then & 
				(x) & if $\varphi \Hderiv \psi$, then $\varphi[t/x] \Hderiv \psi[t/x]$ \\
				& $\varphi \Hderiv \psi \wedge \gamma$; & 
				& ($t$ free for $x$ in $\varphi$ and $\psi$); \\[1ex]
				(iv) & if $\varphi \Hderiv \psi$ and $\psi \Hderiv \gamma$, then $\varphi \Hderiv \gamma$; & 
				(xi) & if $\varphi[c/x] \Hderiv \psi[c/x]$, then $\varphi \Hderiv \psi$ \\
				(v) & if $\varphi \Hderiv \psi$, then $\Diamond\varphi \Hderiv \Diamond\psi$; & 
				& ($c$ not in $\varphi$ nor $\psi$). \\[1ex]
				(vi) & $\Diamond\Diamond\varphi \Hderiv \Diamond\varphi$; & & \\[1ex]
				(vii) & $\Diamond\forall x \varphi \Hderiv \forall x \Diamond\varphi$; & & \\
			\end{tabular}
		\end{center}
		
		\smallskip\noindent
		Axiom (vi) is the \emph{transitivity axiom}. Axiom (vii) is the \emph{quantifier--modality axiom}; its converse $\forall x\Diamond\varphi \Hderiv \Diamond\forall x\varphi$ is not an axiom and is not valid in general.
	\end{definition}
	
	\begin{remark}
		While the base judgments $\varphi \vdash \psi$ of $\QRC$ syntactically resemble sequents, the original calculus is fundamentally a Hilbert-style system. Its rules combine previously established theorems through macroscopic steps such as transitivity and monotonicity. The Hilbert presentation does not directly provide a local root-first proof-search procedure or an analytic subformula discipline. The calculus $\GQRC$ introduced in Section~\ref{sec:calculus} supplies a Gentzen-style presentation and a standard route to cut elimination, the subformula property, and proof-size bounds.
	\end{remark}
	
	\subsection{Kripke semantics\\}
	\label{subsec:kripke}
	
	Our definition of $\QRC$-models streamlines the ``adequate models'' introduced by Borges and Joosten \cite[p.~18--20]{BorgesJoosten2020}. Specifically, we evaluate constants globally rather than by world-indexed interpretations subject to a concordance condition. This yields a simpler semantic presentation, although the two semantics validate exactly the same sequents. This semantic equivalence follows \textit{a posteriori} from our main results: we establish that our sequent calculus is sound and complete with respect to these globally rigid models, and that it is syntactically equivalent to the original Hilbert calculus, which is itself sound and complete with respect to the original adequate models. Thus, our semantics streamline the original models for $\QRC$, as well as the corresponding soundness and completeness proofs.
	
	\begin{definition}[Models for $\QRC$]
		\label{def:model}
		A \emph{$\QRC$-model} is a tuple $\mathcal{M}=(W,R,D,I)$ where:
		\begin{itemize}
			\item $W$ is a non-empty set of \emph{worlds};
			\item $R \subseteq W\times W$ is a \emph{transitive} accessibility relation;
			\item $D: W \to \mathcal{P}^+(\mathcal{U})$ is a \emph{domain function} assigning to each world a non-empty set of \emph{individuals}, where $\mathcal{U}$ is a fixed countably infinite universe of individuals, satisfying:
			\begin{itemize}
				\item \emph{Expanding domains}: if $wRv$ then $D(w) \subseteq D(v)$;
				\item \emph{Constant rigidity}: for each individual constant $c$, $I(c)\in D(w)$ for \emph{every} world $w\in W$. Constants denote inside the local domain at every world.
			\end{itemize}
			\item $I$ assigns to each $n$-ary relation symbol $S$ and world $w$ a set $I_w(S) \subseteq D(w)^n$; and to each individual constant $c$ a \emph{world-independent} element $I(c)\in\bigcap_{w\in W}D(w)$.
		\end{itemize}
		
		A \emph{$w$-assignment} is a function $g_w:\mathsf{Var}\to D(w)$, mapping every variable into the local domain of $w$. Terms are evaluated at $w$ by $\overline{g_w}(x) = g_w(x)$ for variables $x$ and $\overline{g_w}(c) = I(c)$ for constants $c$; by constant rigidity $I(c)\in D(w)$, so $\overline{g_w}$ maps every term to an element of $D(w)$. Since $wRv$ implies $D(w)\subseteq D(v)$, the assignment $g_w$ is automatically a valid $v$-assignment. We will therefore use $g_w$ directly when evaluating formulas at accessible worlds, avoiding the need for assignment coercion. For an assignment \(g_w\), a variable \(x\), and \(d\in D(w)\), we write \(g_w[d/x]\) for the assignment that sends \(x\) to \(d\) and acts as \(g_w\) on all other variables.
		
		The satisfaction relation $\mathcal{M}, w \models_{g_w} \varphi$ is defined inductively:
		\begin{align*}
			\mathcal{M}, w &\models_{g_w} \top \\
			\mathcal{M}, w &\models_{g_w} S(t_1,\ldots,t_n)
			&&\!\!\!\!\text{iff }
			(\overline{g_w}(t_1),\ldots,\overline{g_w}(t_n))\in I_w(S) \\
			\mathcal{M}, w &\models_{g_w} \varphi\wedge\psi
			&&\!\!\!\!\text{iff } \mathcal{M},w\models_{g_w}\varphi
			\text{ and } \mathcal{M},w\models_{g_w}\psi \\
			\mathcal{M}, w &\models_{g_w} \Diamond\varphi
			&&\!\!\!\!\text{iff } \exists v\,(wRv \text{ and }
			\mathcal{M}, v \models_{g_w} \varphi) \\
			\mathcal{M}, w &\models_{g_w} \forall x\varphi
			&&\!\!\!\!\text{iff } \forall d\in D(w)\,
			\mathcal{M},w\models_{g_w[d/x]}\varphi .
		\end{align*}
		
		A sequent $\Gamma\seq C$ is \emph{valid in $\mathcal{M}$} if for all $w\in W$ and all $w$-assignments $g_w$: if $\mathcal{M},w\models_{g_w}\varphi$ for all $\varphi\in\Gamma$, then $\mathcal{M},w\models_{g_w}C$.
	\end{definition}
	
	\begin{lemma}[Coincidence, semantic substitution, and $\alpha$-invariance]
		\label{lem:sem-subst}
		\begin{enumerate}
			\item[(a)] If assignments $g,g'$ at $w$ agree on $\fv(A)$, then
			\[
			\mathcal M,w\models_g A
			\quad\Longleftrightarrow\quad
			\mathcal M,w\models_{g'}A.
			\]
			\item[(b)] For every term $t$,
			\[
			\mathcal M,w\models_g A[t/x]
			\quad\Longleftrightarrow\quad
			\mathcal M,w\models_{g[\overline g(t)/x]}A,
			\]
			where $A[t/x]$ is the capture-avoiding substitution of Definition~\ref{subsdef}.
			\item[(c)] If $A\equiv_\alpha A'$, then
			\[
			\mathcal M,w\models_g A
			\quad\Longleftrightarrow\quad
			\mathcal M,w\models_g A'.
			\]
		\end{enumerate}
	\end{lemma}
	\begin{proof}
		We first prove (a) and (b) simultaneously by strong induction on formula complexity. The cases of $\top$ and conjunction are immediate. In the atomic case, coincidence is immediate, while substitution follows from the elementary term identity
		\[
		\overline g\bigl(s[t/x]\bigr)=\overline{g[\overline g(t)/x]}(s),
		\]
		verified directly according as $s=x$, $s$ is a variable different from $x$, or $s$ is a constant. For $A=\Diamond B$, both assertions follow from the induction hypothesis at each accessible world. If $wRv$, expanding domains keep the relevant assignments admissible at $v$; the value of a variable term under the same assignment function is unchanged when that function is regarded as a $v$-assignment, and the value of a constant term is unchanged by rigidity.
		
		Let $A=\forall zB$. Coincidence follows because, for every $d\in D(w)$, the assignments $g[d/z]$ and $g'[d/z]$ agree on $\fv(B)$. For substitution, if $z=x$, then $(\forall xB)[t/x]=\forall xB$, and coincidence gives
		\[
		\mathcal M,w\models_g\forall xB
		\quad\Longleftrightarrow\quad
		\mathcal M,w\models_{g[\overline g(t)/x]}\forall xB,
		\]
		since $x\notin\fv(\forall xB)$. Now assume $z\neq x$ and put $a=\overline g(t)$. If $z\notin\fv(t)$, then $(\forall zB)[t/x]=\forall z(B[t/x])$ and, for every $d\in D(w)$, $\overline{g[d/z]}(t)=a$. The induction hypothesis and
		\[
		g[d/z][a/x]=g[a/x][d/z]
		\]
		therefore yield the result after quantifying over $d$.
		
		Finally, suppose $z\neq x$ and $z\in\fv(t)$, and let $y$ be the fresh variable chosen in Definition~\ref{subsdef}. For every $d\in D(w)$, two applications of the substitution induction hypothesis give
		\begin{align*}
			\mathcal M,w\models_{g[d/y]}B[y/z][t/x]
			&\Longleftrightarrow \mathcal M,w\models_{g[d/y][a/x]}B[y/z] \\
			&\Longleftrightarrow \mathcal M,w\models_{g[d/y][a/x][d/z]}B \\
			&\Longleftrightarrow \mathcal M,w\models_{g[a/x][d/z]}B.
		\end{align*}
		Here $y\notin\mathrm{Var}(t)$ gives $\overline{g[d/y]}(t)=a$, and $y\neq x$ gives $g[d/y][a/x](y)=d$; the last equivalence is coincidence because $y\notin\mathrm{Var}(B)$. Quantifying over $d$ proves (a) and (b).
		
		For (c), proceed by induction on the common formula structure. The atomic, conjunction, and modal cases are immediate from the induction hypothesis. Suppose
		\[
		A=\forall xB,\qquad A'=\forall yB',
		\]
		and choose $z$ fresh from both formulas. By Definition~\ref{alphadef},
		$B[z/x]\equiv_\alpha B'[z/y]$. For every assignment $g$ and every
		$d\in D(w)$, clauses~(a) and~(b), together with the freshness of $z$,
		give
		\[
		\mathcal M,w\models_{g[d/x]}B
		\quad\Longleftrightarrow\quad
		\mathcal M,w\models_{g[d/z]}B[z/x],
		\]
		and similarly for $B'$ and $y$. The induction hypothesis applied to
		$B[z/x]\equiv_\alpha B'[z/y]$ therefore yields
		\[
		\mathcal M,w\models_{g[d/x]}B
		\quad\Longleftrightarrow\quad
		\mathcal M,w\models_{g[d/y]}B'.
		\]
		Quantifying over $d\in D(w)$ proves
		$\mathcal M,w\models_g\forall xB$ iff
		$\mathcal M,w\models_g\forall yB'$.
	\end{proof}

	\section{The Sequent Calculus \texorpdfstring{$\GQRC$}{G3QRC1}}
	\label{sec:calculus}
	
	\subsection{Sequents\\}
	
	A \emph{sequent} is an expression $\Gamma \seq C$ where $\Gamma$ is a finite multiset of strictly positive formulas (the \emph{antecedent}) and $C$ is a single strictly positive formula (the \emph{succedent}). We write $\Gamma, \Delta$ for the union of multisets and $\Gamma, \varphi$ for adding a single formula $\varphi$ to $\Gamma$. A derivation is a finite tree of sequents built using the rules of the calculus; the \textit{leaves} are axioms and the final sequent at the bottom is called the \emph{root}. A \emph{branch} is a path from the root to a leaf. The \emph{height} $h(\pi)$ of a derivation $\pi$ is the maximum number of sequents along any branch. 
	
	The intended semantic reading of $\Gamma\seq C$ is the $\QRC$-statement $\bigwedge\Gamma \Hderiv C$. Whenever this translation is used, fix an enumeration $A_1,\ldots,A_n$ of the finite multiset $\Gamma$ and let $\bigwedge\Gamma$ denote a fixed parenthesization of $A_1\wedge\cdots\wedge A_n$, with $\bigwedge\emptyset:=\top$. Different enumerations and parenthesizations are interderivable in $\QRC$ by axioms~(ii)--(iv), so the choice is immaterial. The label \textsf{G3} follows the terminology of Troelstra\--Schwichtenberg \cite[p.~77]{TroelstraS2000}: the calculus absorbs all structural rules into the logical rules, so weakening and contraction are admissible.
	
	\subsection{Inference rules of \texorpdfstring{$\GQRC$}{G3QRC1}\\}
	\label{subsec:rules}
	
	The propositional and first-order quantifier rules follow the G3-style presentation of \cite{TroelstraS2000}; the modal rules are new.
	
	\begin{definition}[$\GQRC$]
		The calculus $\GQRC$ has the following axioms and rules.
		
		\medskip
		\noindent\textbf{Identity and verum:}
		\[
		\infer[\rId]{\Gamma, P \seq P}{}
		\qquad
		\infer[\rTopR]{\Gamma \seq \top}{}
		\]
		\textit{In $\rId$, $P$ ranges over \textbf{atomic} formulas only, which are just the relation symbols applied to terms. The admissibility of general identity $\Gamma,A\seq A$ for all formulas $A$ is proved in Lemma~\ref{lem:gen-id}.}
		
		\medskip
		\noindent\textbf{Conjunction:}
		\[
		\infer[\rAndL]{\Gamma, A\wedge B \seq C}{\Gamma, A, B \seq C}
		\qquad
		\infer[\rAndR]{\Gamma \seq A\wedge B}{\Gamma\seq A & \Gamma\seq B}
		\]
		
		\medskip
		\noindent\textbf{Universal quantifier:}
		\[
		\infer[\rAllL]{\Gamma, \forall x\varphi \seq C}
		{\Gamma, \forall x\varphi,\, \varphi[t/x] \seq C}
		\qquad
		\infer[\rAllR]{\Gamma \seq \forall x\varphi}{\Gamma \seq \varphi[y/x]}
		\]
		\textit{Condition on $\rAllR$: the eigenvariable $y$ must not occur free in the conclusion, i.e.\ $y\notin\fv(\Gamma)\cup\fv(\forall x\varphi)$.}
		
		\medskip
		\noindent\textbf{Modal rules:}
		\[
		\infer[\rDiamK]{\Gamma,\Diamond\varphi \seq \Diamond\psi}
		{\varphi \seq \psi}
		\qquad
		\infer[\rDiamFour]{\Gamma,\Diamond\varphi \seq \Diamond\psi}
		{\varphi \seq \Diamond\psi}
		\]
		\textit{In both modal rules the antecedent of the \emph{premise} consists of the \textbf{single} formula $\varphi$. The context $\Gamma$ is introduced entirely in the conclusion.}
		
		\medskip
		\noindent\textbf{Cut:}
		\[
		\infer[\rCut]{\Gamma, \Delta \seq C}
		{\Gamma \seq A & \Delta, A \seq C}
		\]
	\end{definition}
	
	\begin{remark}[Design rationale]
		\label{rem:design}
		\leavevmode
		\begin{enumerate}
			\item[(a)] The $\rAllL$ rule retains $\forall x\varphi$ in the antecedent after instantiation (G3-style), making contraction admissible without a separate rule.
			\item[(b)] The rule $\rDiamK$ encodes monotonicity of $\Diamond$ (Hilbert rule~(v)) directly.
			\item[(c)] The rule $\rDiamFour$ encodes the combination of monotonicity with the transitivity axiom (vi): given $\varphi\Hderiv\Diamond\psi$, one derives $\Diamond\varphi\Hderiv\Diamond\Diamond\psi$ by rule~(v), then applies axiom~(vi) to get $\Diamond\varphi\Hderiv\Diamond\psi$. Encoding this combination as a single primitive rule is exactly what allows the cut-elimination measure to decrease in the principal $\Diamond$-cases (Theorem~\ref{thm:cut}).
			\item[(d)] The context-erasure makes the commutative cut cases for modal rules entirely trivial: a cut formula parametric in the conclusion of a modal rule simply vanishes in the premise, and the modal rule is reapplied with the updated context.
		\end{enumerate}
	\end{remark}
	
	Each of the properties claimed in \ref{rem:design} will be justified.
	
	\section{Structural Properties}
	\label{sec:structural}
	
	Throughout this section, we write $\vdash_n \Gamma \seq C$ to denote that the sequent $\Gamma \seq C$ has a cut-free derivation of height at most $n$. Unless stated otherwise, all derivations are assumed to be cut-free.
	
	Our structural lemmas are established by induction on the derivation height $n$ (often combined with a secondary induction on formula complexity). The standard methodology for these proofs is to proceed by case analysis on the very last inference rule applied to derive what is given, and show how it can be transformed into the desired result. To formalize this analysis, we distinguish between formulas based on their syntactic role within that final rule: \textit{principal formulas} are the specific formulas in the conclusion of an inference rule that are newly introduced or directly operated upon by that rule, whereas \textit{parametric formulas} are the background formulas that are carried from the premises to the conclusion structurally untouched. Consequently, when transforming a derivation to prove a lemma, our analysis naturally splits into \emph{principal cases} (where the final rule directly targets the formula our lemma seeks to manipulate) and \emph{parametric cases} (where the target formula lies untouched in the context). By demonstrating that the derivation can be successfully restructured regardless of which rule was applied last, we guarantee the existence of a valid derivation for the desired conclusion.
	
	\begin{lemma}[General Identity]
		\label{lem:gen-id}
		For any formula $A$, the sequent $\Gamma, A \seq A$ is derivable in $\GQRC$.
	\end{lemma}
	
	\begin{proof}
		By induction on the complexity $\cpx{A}$ of $A$. As our induction hypothesis, we assume that for any formula $F$ with $\cpx{F} < \cpx{A}$, the sequent $\Gamma', F \seq F$ is derivable for any context $\Gamma'$.
		
		\smallskip\noindent\emph{Atomic $P$:} $\Gamma, P\seq P$ is the axiom $\rId$.
		
		\smallskip\noindent\emph{$A = \top$:} $\Gamma, \top\seq\top$ is $\rTopR$.
		
		\smallskip\noindent\emph{$A = B\wedge C$:}
		By the induction hypothesis applied to $B$ with the extended context $\Gamma, C$, we obtain $\Gamma, C, B \seq B$. Applying the induction hypothesis to $C$ with context $\Gamma, B$ gives $\Gamma, B, C \seq C$. From these two sequents, $\rAndR$ yields $\Gamma, B, C \seq B \wedge C$. Finally, $\rAndL$ applied to this premise gives $\Gamma, B\wedge C \seq B\wedge C$, as required.
		
		\smallskip\noindent\emph{$A = \forall xB$:} Let $y$ be a fresh variable with $y\notin\fv(\Gamma,\forall xB)$. The formula $B[y/x]$ has strictly smaller complexity than $\forall xB$; hence the induction hypothesis applies with context $\Gamma,\forall xB$, giving $\Gamma, \forall xB, B[y/x] \seq B[y/x]$. Applying $\rAllL$ with term $y$ yields $\Gamma, \forall xB \seq B[y/x]$. Since $y\notin\fv(\Gamma,\forall xB)$, we may apply $\rAllR$ to obtain $\Gamma, \forall xB \seq \forall xB$.
		
		\smallskip\noindent\emph{$A = \Diamond B$:}
		Because $\cpx{B} < \cpx{\Diamond B}$, the induction hypothesis applies with empty context, yielding $B \seq B$. Applying $\rDiamK$ with context $\Gamma$ to this premise gives $\Gamma, \Diamond B \seq \Diamond B$.
	\end{proof}
	
	\begin{lemma}[$\alpha$-Equivalence]
		\label{lem:alphaeq}
		Let $\varphi$ and $\varphi'$ be formulas that are $\alpha$-equivalent. Then for any context $\Gamma$, the sequents $\Gamma, \varphi \seq \varphi'$ and $\Gamma, \varphi' \seq \varphi$ are derivable.
	\end{lemma}
	
	\begin{proof}
		By structural induction on formula complexity. Our induction hypothesis assumes that for any formula $\psi$ strictly smaller than $\varphi$, and any $\psi'$ that is $\alpha$-equivalent to $\psi$, the sequent $\Delta, \psi \seq \psi'$ is derivable for \emph{any} arbitrary context $\Delta$. For atomic formulas and $\top$, $\alpha$-equivalence implies strict syntactic identity, making $\Gamma, \varphi \seq \varphi'$ an immediate instance of $\rId$ or $\rTopR$. For conjunction, assuming $\varphi = A \wedge B$ and $\varphi' = A' \wedge B'$, the induction hypothesis provides derivations of $\Gamma, B, A \seq A'$ and $\Gamma, A, B \seq B'$; applying $\rAndR$ followed by $\rAndL$ yields $\Gamma, A \wedge B \seq A' \wedge B'$. For the modal case $\varphi = \Diamond A$ and $\varphi' = \Diamond A'$, applying the induction hypothesis with an empty context gives $A \seq A'$; applying $\rDiamK$ then introduces the arbitrary context $\Gamma$ directly in the conclusion to yield $\Gamma, \Diamond A \seq \Diamond A'$. Finally, suppose $\varphi = \forall z \gamma$ and $\varphi' = \forall z' \gamma'$. We choose a variable $y$ that does not occur in any formula of $\Gamma$ or in $\varphi,\varphi'$. By the recursive definition of $\alpha$-equivalence (Definition~\ref{alphadef}), $\forall z \gamma \equiv_\alpha \forall z' \gamma'$ implies that instantiating their bodies with this fresh variable yields $\alpha$-equivalent formulas: $\gamma[y/z] \equiv_\alpha \gamma'[y/z']$. Furthermore, since variable substitution does not change the number of logical connectives, these subformulas remain strictly smaller in complexity than $\forall z \gamma$. We may therefore apply the induction hypothesis to these subformulas. Because our hypothesis applies universally to any context, we specialize it to the extended context $\Delta = \Gamma, \forall z \gamma$, which yields the sequent $\Gamma, \forall z \gamma,\ \gamma[y/z] \seq \gamma'[y/z']$. Applying $\rAllL$ with the term $y$ absorbs the instantiated formula, resulting in $\Gamma, \forall z \gamma \seq \gamma'[y/z']$. Because $y$ was chosen to be completely fresh, it safely satisfies the eigenvariable condition $y \notin \fv(\Gamma, \forall z\gamma)$, allowing a final application of $\rAllR$ to conclude $\Gamma, \forall z \gamma \seq \forall z' \gamma'$. The converse derivation $\Gamma, \varphi' \seq \varphi$ is strictly symmetric.
	\end{proof}
	
	\begin{remark}[$\alpha$-conversion convention]
		\label{rem:alpha-classes}
		Henceforth formulas and sequents are identified up to $\alpha$-equivalence. By Lemma~\ref{lem:alpha-subst}, substitution is well defined on $\alpha$-equivalence classes and bound-variable representatives may be chosen away from any prescribed finite set. The quotient presentation is height-preservingly conservative for $\GQRC$. Indeed, fix syntactic representatives of the end-sequent and proceed by induction on derivation height. At each inference, choose premise representatives according to the selected conclusion representatives and the displayed instantiation terms (for $\rCut$, choose a common representative of the cut formula). Lemmas~\ref{lem:alpha-subst} and~\ref{lem:subst-comp} ensure that these belong to the required $\alpha$-equivalence classes, while eigenvariable conditions remain unchanged. Applying the induction hypothesis to the premises reconstructs a derivation of the same height. In particular, General Identity implies that $\alpha$-equivalent formulas are interderivable. The convention is conservative for the Hilbert calculus as well. If $B\Hderiv C$, then $\forall xB\Hderiv\forall xC$ by rules~(ix) with $t=x$ and~(viii). We prove interderivability of $\alpha$-equivalent formulas by structural induction. The atomic, conjunction, and modal cases are immediate. Suppose $\forall xA\equiv_\alpha\forall yB$. By Lemma~\ref{lem:alpha-subst}, choose $A'\equiv_\alpha A$ whose bound variables avoid $y$. The induction hypothesis gives $A\Hderiv A'$, hence $\forall xA\Hderiv\forall xA'$. Rule~(ix) yields $\forall xA'\Hderiv A'[y/x]$. Since $A'[y/x]\equiv_\alpha B$, the induction hypothesis gives $A'[y/x]\Hderiv B$. By transitivity, $\forall xA\Hderiv B$, and since $y\notin\fv(\forall xA)$, rule~(viii) yields $\forall xA\Hderiv\forall yB$. The converse is symmetric. Accordingly, side conditions such as ``free for'' may be checked after standardizing bound variables apart.
	\end{remark}
	
	Consequently, we may always choose bound variables distinct from any finite collection of free variables or substitution terms under consideration.
	
	\begin{lemma}[Variable Renaming]
		\label{lem:rename}
		Let $\rho$ be an injective map on the set of variables. For a formula $\varphi$, let $\rho\varphi$ denote the result of applying $\rho$ to every variable occurrence in $\varphi$, both free and bound. If $\vdash_n \Gamma \seq C$, then $\vdash_n \rho\Gamma \seq \rho C$. Moreover, uniformly renaming a chosen derivation preserves its height exactly.
	\end{lemma}
	
	\begin{proof}
		By induction on the chosen derivation. The base cases (the axioms $\rId$ and $\rTopR$) are immediate, as renaming maps an axiom to an axiom of the same form. For the inductive step, consider the last rule of the derivation. For each rule ($\rAndL$, $\rAndR$, $\rAllL$, $\rAllR$, $\rDiamK$, $\rDiamFour$) we apply the induction hypothesis to every premise, thereby obtaining renamed derivations of the corresponding premises. The renaming $\rho$ commutes with the logical connectives:
		$\rho(A\land B)=\rho A\land\rho B$ and
		$\rho(\forall z\varphi)=\forall(\rho z)(\rho\varphi)$; for term
		substitution,
		\[
		\rho(\varphi[s/z])\equiv_\alpha
		(\rho\varphi)[\rho s/\rho z].
		\]
		The possible discrepancy comes only from the deterministic choice of
		fresh bound-variable names, and the standing $\alpha$-conversion
		convention therefore suffices to reconstruct the same rule instance. For $\rAllR$, the eigenvariable condition is preserved: $u\notin\fv(\Delta)$ iff $\rho u\notin\fv(\rho\Delta)$, which follows directly from $\fv(\rho\Delta)=\rho(\fv(\Delta))$, a consequence of injectivity. The modal rules $\rDiamK$ and $\rDiamFour$ carry no eigenvariable conditions and simply act on the principal formula and context; they therefore commute with renaming without additional side conditions. Reapplying the same rule to the renamed premises then yields the desired conclusion $\rho\Gamma \seq \rho C$, and no branch length changes.
	\end{proof}
	
	\begin{lemma}[Eigenvariable Renaming]
		\label{lem:eigenrename}
		Let $\pi$ be a cut-free derivation, let $V\subseteq\mathsf{Var}$ be finite, and fix in $\pi$ an application
		\[
		\infer[\rAllR]{\Gamma\seq\forall xA}{\pi_0:\Gamma\seq A[y/x]}.
		\]
		Then there exist a variable $z\in\mathsf{Var}\setminus V$ occurring nowhere in $\pi$ and a cut-free derivation $\pi'$ such that $h(\pi')=h(\pi)$, the end-sequent of $\pi'$ is $\alpha$-equivalent to that of $\pi$, and, after choosing $\alpha$-equivalent representatives, the selected inference in $\pi'$ is
		\[
		\infer[\rAllR]{\Gamma\seq\forall xA}{\pi_0':\Gamma\seq A[z/x]},
		\]
		with $h(\pi_0')=h(\pi_0)$. Every sequent strictly below the selected inference is $\alpha$-equivalent to the corresponding sequent of $\pi$. Consequently, for every finite $V\subseteq\mathsf{Var}$, each cut-free derivation can be transformed into a cut-free derivation of the same height and with an $\alpha$-equivalent end-sequent, in which every $\rAllR$ eigenvariable lies outside $V$.
	\end{lemma}
	\begin{proof}
		By the standing $\alpha$-conversion convention, choose representatives of the displayed inference and its surrounding derivation in which $y$ occurs nowhere in $\Gamma$ or $\forall xA$; if $x=y$, first rename the displayed binder. Choose $z\in\mathsf{Var}\setminus V$ absent from both the original derivation and these representatives, and let $\rho$ be the transposition interchanging $y$ and $z$. Since neither $y$ nor $z$ occurs in the displayed conclusion, $\rho$ fixes $\Gamma$, $x$, and $A$. Lemma~\ref{lem:rename} therefore transforms $\pi_0$ into a derivation $\pi_0'$ of the same height ending in $\Gamma\seq A[z/x]$. Reapplying $\rAllR$ gives the same displayed conclusion, and the part of the derivation strictly below this inference is reattached unchanged at the level of $\alpha$-classes. Thus the total height is preserved.
		
		For the final assertion, repeatedly apply the local construction to an $\rAllR$ occurrence whose eigenvariable belongs to $V$, always choosing the replacement absent from the current derivation and outside $V$. Such a step removes the selected occurrence and introduces no eigenvariable from $V$, so the number of offending occurrences strictly decreases. Since the derivation is finite, the process terminates.
	\end{proof}

	\begin{lemma}[Term Substitution]
		\label{lem:subst}
		If $\vdash_n \Gamma \seq C$, then $\vdash_n \Gamma[t/x] \seq C[t/x]$ for any term $t$ free for $x$ in all formulas in $\Gamma$ and in $C$.
	\end{lemma}
	\begin{proof}
		By induction on the derivation height $n$. Put
		$V=\{x\}\cup\mathrm{Var}(t)$. By Lemma~\ref{lem:eigenrename}, we may
		assume, without changing the height, that every $\rAllR$ eigenvariable
		avoids $V$. Coherently $\alpha$-convert the finite derivation so that
		every bound variable occurring in it avoids $V$. This preserves its
		height and ensures that $t$ is free for $x$ in every formula occurring
		in the derivation, so the induction hypothesis applies to every
		premise. Assume the claim for derivations of height strictly less than
		$n$.
		\begin{itemize}
			\item[$\rId$:] $\Gamma[t/x], A[t/x] \seq A[t/x]$ is an instance of $\rId$.
			\item[$\rTopR$:] $\Gamma[t/x] \seq \top$ is $\rTopR$.
			\item[$\rAndL$:] The premise is $\Gamma, A, B \seq C$. By the induction hypothesis, $\Gamma[t/x], \allowbreak A[t/x], B[t/x] \seq C[t/x]$. Since $(A\wedge B)[t/x] = A[t/x]\wedge B[t/x]$, applying $\rAndL$ yields the required conclusion with height at most $n$.
			\item[$\rAndR$:] Apply the induction hypothesis to both premises and reapply $\rAndR$.
			\item[$\rAllR$:] The premise is $\Gamma \seq \varphi[y/z]$, where $y \notin \fv(\Gamma,\forall z\varphi)$. The preliminary choices give $z\neq x$ and $y,z\notin\fv(t)$. The induction hypothesis yields
			\[
			\Gamma[t/x]\seq(\varphi[y/z])[t/x]
			=(\varphi[t/x])[y/z].
			\]
			Moreover $y\notin\fv(\Gamma[t/x],\forall z(\varphi[t/x]))$, because $y$ is free in neither the original conclusion nor $t$. Hence $\rAllR$ gives
			$\Gamma[t/x]\seq\forall z(\varphi[t/x])=(\forall z\varphi)[t/x]$.
			\item[$\rAllL$:] The premise is $\Gamma,\forall z\varphi,\varphi[s/z]\seq C$. The preliminary choices give $z\neq x$ and $z\notin\fv(t)$. The induction hypothesis yields
			\[
			\Gamma[t/x],\forall z(\varphi[t/x]),(\varphi[t/x])[s[t/x]/z]\seq C[t/x],
			\]
			using $(\varphi[s/z])[t/x]=(\varphi[t/x])[s[t/x]/z]$. Reapply $\rAllL$ with term $s[t/x]$.
			\item[$\rDiamK$, $\rDiamFour$:] Apply the induction hypothesis to the context-free premise and reapply the same modal rule with context $\Gamma[t/x]$.
			\qedhere
		\end{itemize}
	\end{proof}
	
	\begin{lemma}[Constant Substitution]
		\label{lem:const-subst}
		If $\vdash_n \Gamma \seq C$, then $\vdash_n \Gamma[t/c] \seq C[t/c]$ for any constant $c$ and any term $t$, where substitution is capture-avoiding as defined after Definition~\ref{subsdef}.
	\end{lemma}
	\begin{proof}
		Use the same height induction as in Lemma~\ref{lem:subst}. By
		Lemma~\ref{lem:eigenrename}, we may assume that every $\rAllR$
		eigenvariable avoids $\mathrm{Var}(t)$. Coherently $\alpha$-convert
		the finite derivation so that every bound variable occurring in it
		avoids $\mathrm{Var}(t)$. This preserves its height and prevents
		capture throughout the induction. The propositional and modal cases
		are identical to those of Lemma~\ref{lem:subst}. In the quantifier cases one has
		\[
		(\forall z\varphi)[t/c]=\forall z(\varphi[t/c]),
		\qquad
		(\varphi[s/z])[t/c]=(\varphi[t/c])[s[t/c]/z].
		\]
		The induction hypothesis followed by the same quantifier rule therefore yields the required derivation, and the eigenvariable condition is preserved because no eigenvariable occurs in $t$. The height does not increase.
	\end{proof}
	
	\begin{lemma}[Height-Preserving $\top$-Strengthening]
		\label{lem:top}
		If $\vdash_n \Gamma, \top \seq C$, then $\vdash_n \Gamma \seq C$.
	\end{lemma}
	
	\begin{proof}
		By induction on $n$. As our induction hypothesis, we assume that any sequent derivable with height strictly less than \(n\) that contains \(\top\) in the antecedent remains derivable after \(\top\) is removed, without increasing the derivation height.
		\begin{itemize}
			\item $C = \top$: Then $\Gamma \seq \top$ is an instance of $\rTopR$.
			\item Last rule $\rId$: The conclusion must be an instance $\Delta, P \seq P$ with $P$ atomic. Hence $C = P$ and $P$ belongs to the multiset $\Gamma, \top$. Because $\top$ is not atomic, $P \neq \top$; therefore $P$ is already a member of $\Gamma$. Write $\Gamma = \Gamma', P$; then $\Gamma \seq P$ is itself an instance of $\rId$ with context $\Gamma'$.
			\item Last rule $\rDiamK$ or $\rDiamFour$: The distinguished occurrence of $\top$ belongs to the side context. Reapply the same modal rule to the unchanged premise, omitting this occurrence from the side context. The height does not increase.
			\item Last rule is $\rAndL$, $\rAndR$, $\rAllL$, or $\rAllR$: In each case, the formula $\top$ occurs only in the context of the premises. Applying the induction hypothesis to every premise removes $\top$, and then the same rule yields the desired sequent with a derivation of height at most $n$.
			\qedhere
		\end{itemize}
	\end{proof}
	
	\begin{lemma}[Height-Preserving Weakening]
		\label{lem:weak}
		If $\vdash_n \Gamma \seq C$, then $\vdash_n \Gamma, \Delta \seq C$ for any finite multiset $\Delta$.
	\end{lemma}
	
	\begin{proof}
		Induction on $n$. As our induction hypothesis, we assume that any sequent derivable with height strictly less than $n$ can be weakened by any finite multiset without increasing the derivation height.
		\begin{itemize}
			\item $\rId$, $\rTopR$: the additional formulas in $\Delta$ are vacuously admitted; $\Gamma,\Delta\seq P$ (atomic $P\in\Gamma$) and $\Gamma,\Delta \seq \top$ are instances of the respective axioms.
			\item $\rDiamK$, $\rDiamFour$: write $\Gamma=\Gamma_0,\Diamond\varphi$, distinguishing the principal occurrence. The premise $\varphi\seq\psi$ (or $\varphi\seq\Diamond\psi$) is unchanged. Reapply the same modal rule with side context $\Gamma_0,\Delta$. The resulting derivation has height at most $n$.
			\item $\rAndL$, $\rAndR$, $\rAllL$, $\rAllR$: apply the induction hypothesis to each premise and reapply the rule. For $\rAllR$, if its eigenvariable $y$ occurs free in $\Delta$, apply the local clause of Lemma~\ref{lem:eigenrename} to that final inference with $V=\fv(\Delta)$, obtaining a replacement eigenvariable $y'\notin\fv(\Delta)$, and then reapply $\rAllR$. The resulting derivation has height at most $n$.
			\qedhere
		\end{itemize}
	\end{proof}
	
	\begin{lemma}[Height-Preserving $\wedge L$-Invertibility]
		\label{lem:andl-inv}
		If $\vdash_n \Gamma, A\wedge B \seq C$, then $\vdash_n \Gamma, A, B \seq C$.
	\end{lemma}
	
	\begin{proof}
		Induction on $n$. As our induction hypothesis, we assume that for any derivation of height strictly less than \(n\) ending in \(\Gamma', A\land B \seq C'\), the sequent \(\Gamma', A, B \seq C'\) is derivable with the same height bound.
		\begin{itemize}
			\item If the last rule is $\rId$ (so $C = P$ is an atomic formula and $P\in\Gamma,A\wedge B$): since $A\wedge B$ is not atomic, $P\in\Gamma$. The sequent $\Gamma, A, B \seq P$ is an instance of $\rId$ at height 1. 
			\item If the last rule is $\rTopR$ (so $C = \top$): $\Gamma, A, B\seq\top$ by $\rTopR$ at height 1. 
			\item If the last rule is $\rAndL$ with $A\wedge B$ principal: the premise is $\Gamma, A, B \seq C$ at height at most $n-1$.
			\item If the last rule is $\rDiamK$ or $\rDiamFour$: the distinguished occurrence of $A\wedge B$ belongs to the side context. Reapply the same modal rule to the unchanged premise, replacing this occurrence in the side context by $A,B$. The height does not increase.
			\item Otherwise (last rule is $\rAndL$ with a principal formula different from $A\land B$, or $\rAndR$, $\rAllL$, $\rAllR$): $A\wedge B$ is parametric; apply the induction hypothesis to each premise and reapply the rule. For $\rAllR$, the eigenvariable condition remains valid because the set of free variables is unchanged: $\fv(A, B) = \fv(A \wedge B)$. The height does not increase.
			\qedhere
		\end{itemize}
	\end{proof}
	
	\begin{theorem}[Admissibility of Contraction]
		\label{thm:contr}
		If $\vdash \Gamma, A, A \seq C$, then $\vdash \Gamma, A \seq C$.
	\end{theorem}
	
	\begin{proof}
		We argue by double induction: the outer induction on the complexity $\cpx{A}$ of the contracted formula, and the inner induction on the height $n$ of the given cut-free derivation of $\Gamma, A, A \seq C$. As our outer induction hypothesis, we assume that for every formula $F$ with $\cpx{F} < \cpx{A}$, any sequent $\Gamma', F, F \seq C'$ can be contracted to $\Gamma', F \seq C'$. As our inner induction hypothesis, we assume this contraction property holds for the formula $A$ itself, provided the derivation has height strictly smaller than $n$.
		
		We distinguish whether the last rule of the given derivation acts on one of the two copies of $A$ (principal case) or leaves both copies untouched in the context (parametric case).
		
		\smallskip\noindent
		\emph{Axioms.}
		If the derivation is the axiom $\rId$, then $C = P$ is atomic and $P$ belongs to the multiset $\Gamma, A, A$. When $P = A$, the sequent $\Gamma, A \seq A$ is an instance of the $\rId$ axiom itself. If $P \in \Gamma$, then $\Gamma, A \seq P$ is again the $\rId$ axiom. The $\rTopR$ axiom yields $\Gamma, A \seq \top$ trivially.
		
		\smallskip\noindent
		\emph{Principal case for conjunction.}
		Suppose $A = B \land D$ and the last rule is $\rAndL$ applied to one copy of $B \land D$, say the first. The premise is
		\[
		\vdash_{n-1} \Gamma,\; B\land D,\; B,\; D \seq C .
		\]
		By height-preserving invertibility of $\rAndL$ (Lemma~\ref{lem:andl-inv}), we replace the remaining $B\land D$ in the antecedent by $B, D$, obtaining a derivation of
		\[
		\vdash_{n-1} \Gamma,\; B,\; D,\; B,\; D \seq C .
		\]
		Since $\cpx{B}, \cpx{D} < \cpx{B\land D}$, the outer induction hypothesis applies. Contracting first the two copies of $B$ and then the two copies of $D$ yields $\vdash \Gamma, B, D \seq C$. A final application of $\rAndL$ gives $\Gamma, B\land D \seq C$, as required.
		
		\smallskip\noindent
		\emph{Principal case for universal quantifier.}
		Assume $A = \forall z\,\varphi$ and the last rule is $\rAllL$ that instantiates one copy of $\forall z\,\varphi$ with the term $t$. The premise is
		\[
		\vdash_{n-1} \Gamma,\; \forall z\varphi,\; \forall z\varphi,\; \varphi[t/z] \seq C .
		\]
		Here the inner induction hypothesis (same formula, height $n-1$) directly contracts the two copies of $\forall z\varphi$, giving
		\[
		\vdash \Gamma,\; \forall z\varphi,\; \varphi[t/z] \seq C.
		\]
		Applying $\rAllL$ again yields the contracted sequent $\Gamma, \forall z\varphi \seq C$.
		
		\smallskip\noindent
		\emph{Principal case for diamond.}
		Let $A = \Diamond\varphi$ and suppose the last rule is $\rDiamK$ or $\rDiamFour$ that introduces one copy of $\Diamond\varphi$ as its principal formula. Both rules have the form
		\[
		\infer{\Delta, \Diamond\varphi \seq \Diamond\psi}{\varphi \seq \psi}
		\qquad\text{or}\qquad
		\infer{\Delta, \Diamond\varphi \seq \Diamond\psi}{\varphi \seq \Diamond\psi},
		\]
		where $\Delta$ collects the parametric context. In the conclusion $\Gamma, \Diamond\varphi, \Diamond\varphi \seq \Diamond\psi$ one occurrence is principal and the other belongs to $\Delta$. We may therefore apply the \emph{same} rule again, but now with context $\Gamma$ instead of $\Gamma, \Diamond\varphi$. The premise does not change, so no induction hypothesis is needed; the contracted sequent $\Gamma, \Diamond\varphi \seq \Diamond\psi$ is obtained immediately.
		
		\smallskip\noindent
		\emph{Parametric cases.}
		If the last rule is a modal rule $\rDiamK$ or $\rDiamFour$ and neither distinguished occurrence of $A$ is principal, both belong to the context of the conclusion and are erased in the premise. We simply apply the same modal rule to the unchanged premise but with one of the two distinguished occurrences of $A$ deleted from its side context, directly obtaining the contracted conclusion.
		
		If the last rule is any other non-modal rule ($\rAndL$ with a different principal formula, $\rAndR$, $\rAllL$, $\rAllR$), then both copies of $A$ appear parametrically in every premise. The inner induction hypothesis applies to each premise (same $A$, strictly smaller height) and we afterward reapply the same rule. For the rule $\rAllR$ one simply notes that the eigenvariable condition remains satisfied: the fresh eigenvariable $y$ does not belong to $\fv(\Gamma, A, A)$, and since $\fv(\Gamma, A, A) = \fv(\Gamma, A)$, it is also not free in the contracted context $\Gamma, A$.
		
		All possibilities for the last rule have been examined, so the induction is complete.
	\end{proof}

	\section{Cut-Elimination}
	\label{sec:cut}
	
	\begin{theorem}[Cut Admissibility]
		\label{thm:cut}
		If $\Gamma\seq A$ and $\Delta,A\seq C$ have cut-free derivations in $\GQRC$, then so does $\Gamma,\Delta\seq C$. Hence every derivable sequent in $\GQRC$ has a cut-free derivation.
	\end{theorem}
	
	\begin{proof}
		We first reduce to the case of a \emph{topmost cut}: a cut whose two premises are themselves cut-free. Any derivation containing at least one cut contains a topmost cut. By an outer induction on the total number of cuts in a derivation, it therefore suffices to eliminate a single topmost cut. With the premises $\Gamma\seq A$ and $\Delta,A\seq C$ cut-free, every lemma from Section~\ref{sec:structural} applies to them and to their sub-derivations, since those lemmas were proved for cut-free derivations.
		
		We now fix such a topmost cut, with cut-free left premise of height $h_L$ and cut-free right premise of height $h_R$. We prove that this cut can be eliminated by lexicographic induction on $\langle \cpx{A},\, h_L + h_R \rangle$. The induction hypothesis is:
		\begin{quote}
			For any cut with cut formula $B$ and premises heights $h'_L, h'_R$, if $\langle \cpx{B},\, h'_L + h'_R \rangle < \langle \cpx{A},\, h_L + h_R \rangle$, then a cut-free derivation of the corresponding conclusion exists.
		\end{quote}
		Here the lexicographic order on pairs means that either $\cpx{B}<\cpx{A}$, or $\cpx{B}=\cpx{A}$ and $h'_L+h'_R < h_L+h_R$. Because the original premises are cut-free, and every recursive call preserves cut-freeness as the hypothesis returns a cut-free derivation, all lemmas of Section~\ref{sec:structural} are applicable throughout the argument.
		
		We now analyse the possible last rules of the left and right premises. Throughout, ``the cut formula'' means the distinguished occurrence introduced into the right premise of the cut; all multiset decompositions and principal/parametric distinctions are occurrence-sensitive. As in the structural lemmas, we distinguish whether that occurrence is \emph{principal} in a premise or \emph{parametric}. The induction hypothesis is invoked whenever the measure decreases; when the measure does not decrease, the case is reduced to a trivial transformation (as in the axiom cases) or to an application of a modal rule that resets the context.
		
		\subsection*{Group 1: Axiom Cases}
		
		\paragraph{Case 1a: Left premise is $\rId$.}
		The left derivation is $\Gamma_0, A \seq A$ (height 1), so $\Gamma = \Gamma_0, A$. The cut conclusion is $\Gamma_0, A, \Delta \seq C$. This is precisely Lemma~\ref{lem:weak} applied to the right premise $\Delta, A \seq C$ (adding $\Gamma_0$): $\vdash \Gamma_0, \Delta, A \seq C$, i.e., $\vdash \Gamma_0, A, \Delta \seq C$.
		
		\paragraph{Case 1b: Right premise is $\rId$.}
		The right premise $\Delta, A \seq C$ is an $\rId$ axiom instance. Since $\rId$ is restricted to atomic formulas, $C$ must be atomic and $C \in \Delta \cup \{A\}$. We split into two sub-cases.
		
		\smallskip
		\noindent\textit{Sub-case (i): $C \in \Delta$ (cut formula $A$ parametric in $\rId$).} Write $\Delta = \Delta', C$. The cut conclusion is $\Gamma, \Delta', C \seq C$, which is an instance of $\rId$.
		
		\smallskip
		\noindent\textit{Sub-case (ii): $C = A$ (cut formula $A$ is principal in $\rId$, $C \notin \Delta$).} Then $A$ is atomic and the right premise is $\Delta, A \seq A$. The cut conclusion is $\Gamma, \Delta \seq A = \Gamma, \Delta \seq C$. Apply Weakening (Lemma~\ref{lem:weak}) to the left premise $\Gamma \seq A$ (adding $\Delta$): $\vdash \Gamma, \Delta \seq A$.
		
		\paragraph{Case 1c: Cut formula is $\top$.} The left premise is unused, since its end-sequent $\Gamma\seq\top$ is independently an instance of $\rTopR$. The right premise is $\Delta, \top \seq C$. By Lemma~\ref{lem:top} ($\top$-Strengthening) on the right premise, $\vdash \Delta \seq C$ without increasing the height. Apply Lemma~\ref{lem:weak}: $\vdash \Gamma, \Delta \seq C$.
		
		\paragraph{Case 1d: Right premise is $\rTopR$.} Right premise: $\Delta, A \seq \top$ (by $\rTopR$, so $C = \top$). The cut conclusion $\Gamma, \Delta \seq \top$ is immediately an instance of $\rTopR$. No use of the left premise or induction hypothesis is needed.
		
		\subsection*{Group 2: Commutative (Permutation) Cases} In this group of cases, the cut formula $A$ is parametric with respect to the last inference rule applied in the indicated premise. Because this preceding rule operates strictly on formulas in the background context and leaves the cut formula undisturbed, its application is structurally independent of the cut. 
		
		This independence allows us to permute the order of inferences: we can push the $\rCut$ upward, applying it before the background rule rather than after. By performing the cut one step higher in the derivation, the sum of the premise heights $h_L + h_R$ strictly decreases while the cut formula complexity $\cpx{A}$ remains unchanged. This strict reduction in our lexicographic measure allows us to safely invoke the induction hypothesis, which immediately yields a cut-free derivation of the upper branch. The permuted rule is then reapplied at the bottom to yield the final conclusion.
		
		In a commutative case the cut formula $A$ is not principal in the premise under discussion. In the left premise $\Gamma\seq A$, the cut formula is literally the succedent. Hence, if the last rule were $\rAndR$, $\rAllR$, $\rDiamK$, or $\rDiamFour$, that rule would introduce the succedent and would make $A$ principal. Apart from the axiom cases already treated, the last rule on the left must therefore be $\rAndL$ or $\rAllL$, the only rules that carry the succedent parametrically. In the right premise $\Delta,A\seq C$, by contrast, the cut formula lies in the antecedent, so a right rule may leave it parametric, and a modal rule may erase it as part of its context; all such possibilities are treated below.
		\smallskip
		
		\paragraph{Case 2a: Left premise ends in $\rAndL$ or $\rAllL$, $A$ parametric.}
		\begin{itemize}
			\item \emph{$\rAndL$:}
			Left: $\dfrac{\Gamma', D, E \seq A}{\Gamma', D\wedge E \seq A}$ so $\Gamma = \Gamma', D\wedge E$, and the premise has height $h_L-1$. We form the cut
			\[
			\infer[\rCut]{\Gamma', D, E, \Delta \seq C}{
				\Gamma', D, E \seq A & \Delta, A \seq C}
			\]
			whose measure $\langle\cpx{A}, (h_L-1)+h_R\rangle$ is strictly smaller. By the induction hypothesis, we obtain a cut-free derivation of $\Gamma', D, E, \Delta \seq C$. Applying $\rAndL$ yields $\Gamma', D\wedge E, \Delta \seq C$, i.e.\ $\Gamma,\Delta \seq C$.
			\item \emph{$\rAllL$:} Identical logic; push Cut into the premise and reapply $\rAllL$.
		\end{itemize}
		
		\paragraph{Case 2b: Right premise ends in $\rAndL$, $A$ parametric.}
		Right: $\dfrac{\Delta', D, E \seq C}{\Delta', D\wedge E \seq C}$. The right premise of the cut is $\Delta, A \seq C$, which is exactly the conclusion of this $\rAndL$ rule; hence, as multisets with distinguished occurrences, $\Delta,A=\Delta',D\wedge E$. Since the distinguished occurrence of $A$ is parametric, it occurs in $\Delta'$. Write $\Delta' = \Delta_0, A$. Then $\Delta = \Delta_0, D\wedge E$, and the premise of the $\rAndL$ rule (height $h_R-1$) is $\Delta_0, A, D, E \seq C$. We form the cut directly:
		\[
		\infer[\rCut]{\Gamma, \Delta_0, D, E \seq C}{
			\Gamma \seq A \quad \Delta_0, A, D, E \seq C}
		\]
		whose measure $\langle\cpx{A},\, h_L + (h_R-1)\rangle$ is strictly smaller. By the induction hypothesis we obtain a cut-free derivation of $\Gamma, \Delta_0, D, E \seq C$. Applying $\rAndL$ yields $\Gamma, \Delta_0, D\wedge E \seq C$, i.e.\ $\Gamma, \Delta \seq C$.
		
		\paragraph{Case 2c: Right premise ends in $\rAndR$, $A$ parametric.} Right: $\dfrac{\Delta, A \seq D \quad \Delta, A \seq E} {\Delta, A \seq D\wedge E}$. Both premises have height at most $h_R-1$. Form the cuts of $\Gamma \seq A$ with each premise; each cut has the same cut formula $A$ and smaller total height. By the induction hypothesis we obtain cut-free derivations of $\Gamma, \Delta \seq D$ and $\Gamma, \Delta \seq E$. Applying $\rAndR$ yields $\Gamma, \Delta \seq D\wedge E$.
		
		\paragraph{Case 2d: Right premise ends in $\rAllL$, $A$ parametric.} Right: $\dfrac{\Delta', A, \forall z\varphi, \varphi[t/z] \seq C} {\Delta', A, \forall z\varphi \seq C}$, where the premise has height $h_R-1$. Form the cut of $\Gamma \seq A$ with this premise; the cut formula is $A$ and the total height decreases. By the induction hypothesis we obtain a cut-free derivation of $\Gamma, \Delta', \forall z\varphi, \varphi[t/z] \seq C$. Applying $\rAllL$ yields $\Gamma, \Delta', \forall z\varphi \seq C$, which is exactly the desired conclusion $\Gamma, \Delta \seq C$.
		
		\paragraph{Case 2e: Right premise ends in $\rAllR$, $A$ parametric.} Right: $\dfrac{\Delta, A \seq \psi[y/z]}{\Delta, A \seq \forall z\psi}$ with eigenvariable $y\notin\fv(\Delta, A, \forall z\psi)$. First apply Lemma~\ref{lem:eigenrename} with $V=\fv(\Gamma,\Delta,A,\forall z\psi)$ to replace the eigenvariable by a fresh $y'$; the premise becomes $\Delta, A \seq \psi[y'/z]$ (height $h_R-1$). Form the cut of $\Gamma \seq A$ with this premise; the cut formula is $A$ and the total height decreases. By the induction hypothesis we obtain a cut-free derivation of $\Gamma, \Delta \seq \psi[y'/z]$ with $y'\notin\fv(\Gamma,\Delta)$. Applying $\rAllR$ then yields $\Gamma, \Delta \seq \forall z\psi$, which is the desired conclusion.

		\paragraph{Case 2f: Right premise ends in $\rDiamK$ or $\rDiamFour$, $A$ parametric.}
		Right: $\dfrac{\varphi \seq \psi}{\Delta', A, \Diamond\varphi \seq \Diamond\psi}$ where the distinguished occurrence of $A$ belongs to the context and is erased in the premise. To obtain the required cut conclusion $\Gamma, \Delta', \Diamond\varphi \seq \Diamond\psi$, we simply apply the same modal rule again, but now with context $\Gamma, \Delta'$ instead of $\Delta', A$:
		\[
		\infer[\rDiamK]{\Gamma,\Delta',\Diamond\varphi\seq\Diamond\psi}{\varphi\seq\psi}
		\]
		The $\rDiamFour$ case is identical, using the premise $\varphi\seq\Diamond\psi$. No use of the left derivation or the induction hypothesis is required.
		
		\subsection*{Group 3: Principal Cases}
		The previous group considered the cases where the cut formula $A$ was non-principal in either premise. We now consider the cases where $A$ is principal in both premises. If $A$ is a conjunction, the relevant final rules are $\rAndR$ on the left and $\rAndL$ on the right; if $A$ is universal, they are $\rAllR$ on the left and $\rAllL$ on the right; and if $A$ is diamond-headed, either premise may end in either modal rule. In every case the lexicographic cut measure decreases. Some reductions first recur on the same cut formula with a smaller sum of premise heights, as in the universal case and the cases where the left premise ends in $\rDiamFour$; the remaining recursive cuts use formulas of strictly smaller logical complexity, including substitution instances of proper subformulas in the universal case. When the left premise ends in $\rDiamFour$, the right modal conclusion is reconstructed with empty parametric context without increasing its derivation height.
		
		\paragraph{Case 3a: $A = D\wedge E$.}
		The left premise ends with $\rAndR$:
		\[
		\infer[\rAndR]{\Gamma \seq D \wedge E}{\Gamma \seq D & \Gamma \seq E}
		\]
		and the right premise ends with $\rAndL$:
		\[
		\infer[\rAndL]{\Delta, D\wedge E \seq C}{ \Delta, D, E \seq C}.
		\]
		The sub-derivations $\Gamma \seq D$ and $\Gamma \seq E$ have height $\le h_L-1$; the sub-derivation $\Delta, D, E \seq C$ has height $\le h_R-1$. First we form a cut on $E$:
		\[
		\infer[\rCut]{\Gamma, \Delta, D \seq C}{\Gamma \seq E & \Delta, D, E \seq C}
		\]
		Since $\cpx{E} < \cpx{D\wedge E}$, the lexicographic measure of this new cut is strictly smaller than $\langle\cpx{A}, h_L+h_R\rangle$. By the induction hypothesis we obtain a cut-free derivation of $\Gamma, \Delta, D \seq C$. Now we form a second cut, this time on $D$, between the sub-derivation $\Gamma \seq D$ and the just obtained derivation of $\Gamma, \Delta, D \seq C$:
		\[
		\infer[\rCut]{\Gamma, \Gamma, \Delta \seq C}{\Gamma \seq D & \Gamma, \Delta, D \seq C}
		\]
		Again $\cpx{D} < \cpx{D\wedge E}$, so the measure is smaller, and the induction hypothesis gives a cut-free derivation of $\Gamma, \Gamma, \Delta \seq C$. Applying contraction (Theorem~\ref{thm:contr}) finitely many times removes the duplicate copies of $\Gamma$ and yields $\Gamma, \Delta \seq C$.
		
		\paragraph{Case 3b: $A = \forall z\varphi$.} The left premise ends with $\rAllR$:
		\[
		\infer[\rAllR]{\Gamma \seq \forall z\varphi}{\Gamma \seq \varphi[y/z]}
		\]
		where $y \notin \fv(\Gamma,\forall z\varphi)$ and the sub-derivation $\Gamma \seq \varphi[y/z]$ has height $\le h_L-1$. The right premise ends with $\rAllL$:
		\[
		\infer[\rAllL]{\Delta,\forall z\varphi \seq C}{\Delta,\forall z\varphi,\varphi[t/z] \seq C}
		\]
		whose sub-derivation has height $\le h_R-1$. We first apply the induction hypothesis to a new cut on the same formula $\forall z\varphi$ but with a strictly smaller total height. Take as premises the left derivation $\Gamma \seq \forall z\varphi$ (height $h_L$) and the premise of the right $\rAllL$ rule, namely $\Delta,\forall z\varphi,\varphi[t/z] \seq C$ (height $\le h_R-1$). The resulting cut
		\[
		\infer[\rCut]{\Gamma,\Delta,\varphi[t/z] \seq C}{\Gamma \seq \forall z\varphi & \Delta,\forall z\varphi,\varphi[t/z] \seq C}
		\]
		has cut formula $\forall z\varphi$ and total height at most $h_L + (h_R-1) < h_L + h_R$. Since the complexity of the cut formula is unchanged while the sum of heights strictly decreases, the lexicographic measure becomes smaller; the induction hypothesis therefore gives a cut-free derivation of $\Gamma,\Delta,\varphi[t/z] \seq C$. 
		
		It remains to obtain a derivation of $\Gamma \seq \varphi[t/z]$. Apply the local clause of Lemma~\ref{lem:eigenrename} to the final $\rAllR$ inference, taking $V$ to be the finite set of all variables occurring in $\Gamma$, $\varphi$, and $t$, together with $z$. We thereby obtain, without changing the height, a premise derivation of $\Gamma\seq\varphi[y'/z]$ with $y'\notin V$. Before applying term substitution, use the standing
		$\alpha$-conversion convention to standardize every bound variable of
		$\varphi$ apart from $\{z,y'\}\cup\mathrm{Var}(t)$, as required by
		Lemma~\ref{lem:subst-comp}. Then $t$ is free for $y'$ in
		$\varphi[y'/z]$, so Lemma~\ref{lem:subst} applies with $[t/y']$. Since $y'\notin\fv(\Gamma)$, the antecedent remains $\Gamma$, and Lemma~\ref{lem:subst-comp}, with the stated freshness choices, gives
		\[
		(\varphi[y'/z])[t/y']\equiv_\alpha\varphi[t/z].
		\]
		Hence $\Gamma\seq\varphi[t/z]$ has a cut-free derivation of height at most $h_L-1$. We now form a second cut on the strictly smaller formula $\varphi[t/z]$:
		\[
		\infer[\rCut]{\Gamma,\Gamma,\Delta \seq C}{\Gamma \seq \varphi[t/z] & \Gamma,\Delta,\varphi[t/z] \seq C}
		\]
		Since $\cpx{\varphi[t/z]} = \cpx{\varphi} < \cpx{\forall z\varphi}$, the complexity decreases, so the induction hypothesis applies again and gives a cut-free derivation of $\Gamma,\Gamma,\Delta \seq C$. Finally, contraction (Theorem~\ref{thm:contr}) applied finitely many times removes the duplicate copies of $\Gamma$, yielding the desired $\Gamma,\Delta \seq C$.
		
		\paragraph{Case 3c: $A = \Diamond\varphi$, left ends in $\rDiamK$, right ends in $\rDiamK$.}
		The left premise is
		\[
		\infer[\rDiamK]{\Gamma,\Diamond \gamma \seq \Diamond\varphi}{\gamma \seq \varphi}
		\]
		and the right premise is
		\[
		\infer[\rDiamK]{\Delta,\Diamond\varphi \seq \Diamond\psi}{\varphi \seq \psi}.
		\]
		Both sub-derivations ($\gamma \seq \varphi$ and $\varphi \seq \psi$) have heights $\le h_L-1$ and $\le h_R-1$, respectively. We form a cut on the immediate subformula $\varphi$:
		\[
		\infer[\rCut]{\gamma \seq \psi}{\gamma \seq \varphi & \varphi \seq \psi}
		\]
		Since $\cpx{\varphi} < \cpx{\Diamond\varphi}$, the complexity decreases; the induction hypothesis yields a cut-free derivation of $\gamma \seq \psi$. Applying the rule $\rDiamK$ with the combined context $\Gamma,\Delta$ gives the desired $\Gamma,\Delta,\Diamond \gamma \seq \Diamond\psi$.
		
		\paragraph{Case 3d: $A = \Diamond\varphi$, left ends in $\rDiamK$, right ends in $\rDiamFour$.}
		The left premise is
		\[
		\infer[\rDiamK]{\Gamma,\Diamond \gamma \seq \Diamond\varphi}{\gamma \seq \varphi}
		\]
		and the right premise is
		\[
		\infer[\rDiamFour]{\Delta,\Diamond\varphi \seq \Diamond\psi}{\varphi \seq \Diamond\psi}.
		\]
		Both sub-derivations have heights $\le h_L-1$ and $\le h_R-1$, respectively. We form a cut on the immediate subformula $\varphi$:
		\[
		\infer[\rCut]{\gamma \seq \Diamond\psi}{\gamma \seq \varphi & \varphi \seq \Diamond\psi}
		\]
		Since $\cpx{\varphi} < \cpx{\Diamond\varphi}$, the complexity decreases; the induction hypothesis yields a cut-free derivation of $\gamma \seq \Diamond\psi$. Applying $\rDiamFour$ with the combined context $\Gamma,\Delta$ gives the desired $\Gamma,\Delta,\Diamond \gamma \seq \Diamond\psi$.
		
		\paragraph{Case 3e: $A = \Diamond\varphi$, left ends in $\rDiamFour$, right ends in $\rDiamK$.}
		The left premise is
		\[
		\infer[\rDiamFour]{\Gamma,\Diamond \gamma \seq \Diamond\varphi}{\gamma \seq \Diamond\varphi}
		\]
		with sub-derivation $\gamma \seq \Diamond\varphi$ of height $\le h_L-1$. The right premise is
		\[
		\infer[\rDiamK]{\Delta,\Diamond\varphi \seq \Diamond\psi}{\varphi \seq \psi}
		\]
		with sub-derivation $\varphi \seq \psi$ of height $\le h_R-1$. We cannot directly cut $\gamma \seq \Diamond\varphi$ against $\varphi \seq \psi$ because the cut formulas ($\Diamond\varphi$ vs.\ $\varphi$) differ. Instead, we first apply $\rDiamK$ (with empty context) to the sub-derivation $\varphi \seq \psi$ (height $\le h_R-1$), obtaining a derivation of $\Diamond\varphi \seq \Diamond\psi$ with height $\le h_R$. Now we form a cut on $\Diamond\varphi$:
		\[
		\infer[\rCut]{\gamma \seq \Diamond\psi}{\gamma \seq \Diamond\varphi & \Diamond\varphi \seq \Diamond\psi}
		\]
		The measure of this cut is $\langle\cpx{\Diamond\varphi},\,(h_L-1)+h_R\rangle$, which is strictly smaller than the original $\langle\cpx{\Diamond\varphi},\, h_L+h_R\rangle$ because the left premise height decreases by one. Hence the induction hypothesis gives a cut-free derivation of $\gamma \seq \Diamond\psi$. Finally, applying $\rDiamFour$ with context $\Gamma,\Delta$ yields $\Gamma,\Delta,\Diamond\gamma \seq \Diamond\psi$, as required.
		
		\paragraph{Case 3f: $A = \Diamond\varphi$, left ends in $\rDiamFour$, right ends in $\rDiamFour$.}
		The left premise is
		\[
		\infer[\rDiamFour]{\Gamma,\Diamond \gamma \seq \Diamond\varphi}{\gamma \seq \Diamond\varphi}
		\]
		with sub-derivation $\gamma \seq \Diamond\varphi$ of height $\le h_L-1$. The right premise is
		\[
		\infer[\rDiamFour]{\Delta,\Diamond\varphi \seq \Diamond\psi}{\varphi \seq \Diamond\psi}
		\]
		with sub-derivation $\varphi \seq \Diamond\psi$ of height $\le h_R-1$. We first apply $\rDiamFour$ with empty context to the right sub-derivation $\varphi \seq \Diamond\psi$, obtaining a derivation of $\Diamond\varphi \seq \Diamond\psi$ with height $\le h_R$. Now we cut on $\Diamond\varphi$:
		\[
		\infer[\rCut]{\gamma \seq \Diamond\psi}{\gamma \seq \Diamond\varphi & \Diamond\varphi \seq \Diamond\psi}
		\]
		The measure of this cut is $\langle\cpx{\Diamond\varphi},\, (h_L-1)+h_R\rangle$, which is strictly smaller than the original $\langle\cpx{\Diamond\varphi},\, h_L+h_R\rangle$. Hence the induction hypothesis gives a cut-free derivation of $\gamma \seq \Diamond\psi$. Finally, applying $\rDiamFour$ with context $\Gamma,\Delta$ yields $\Gamma,\Delta,\Diamond\gamma \seq \Diamond\psi$, as required.
		
		\bigskip
		All cases are covered. The lexicographic measure $\langle\cpx{A}, h_L+h_R\rangle$ decreases at every recursive call, so the induction terminates. Cut is admissible.
	\end{proof}

	\section{Soundness and Completeness}
	\label{sec:completeness}
	
	\subsection{Soundness\\}
	
	We interpret a sequent $\Gamma\seq C$ semantically as: for all $\QRC$-models $\mathcal{M}=(W,R,D,I)$, all $w\in W$, and all $w$-assignments $g_w:\mathsf{Var}\to D(w)$, $(\mathcal{M},w\models_{g_w}\varphi\text{ for all }\varphi\in\Gamma)\;\Rightarrow\; \mathcal{M},w\models_{g_w}C$.
	
	\begin{theorem}[Soundness]
		\label{thm:sound}
		If $\vdash \Gamma\seq C$ in $\GQRC$, then $\Gamma\seq C$ is valid in every $\QRC$-model.
	\end{theorem}
	
	\begin{proof}
		By Theorem~\ref{thm:cut} (Cut Admissibility), every derivable sequent has a cut-free derivation. Therefore it suffices to prove soundness for cut-free derivations; we proceed by induction on the height of a cut-free derivation, where the induction hypothesis is that every sequent with a cut-free derivation of height strictly less than $n$ is valid in all $\QRC$-models. We now examine the last rule used in a given derivation.
		
		\smallskip
		\noindent\textbf{Axiom $\rId$:} The sequent is $\Gamma, P \seq P$ with $P$ atomic. Assume that $\mathcal{M},w \models_{g_w} \bigwedge\Gamma \wedge P$. Then in particular $\mathcal{M},w \models_{g_w} P$, so the succedent holds. Hence the sequent is valid.
		
		\smallskip
		\noindent\textbf{Axiom $\rTopR$:} The sequent $\Gamma \seq \top$ is valid because $\mathcal{M},w \models_{g_w} \top$ holds for every world $w$ and every assignment $g_w$.
		
		\smallskip
		\noindent\textbf{Rule $\rAndL$:} The premise is $\Gamma, A, B \seq C$, derived with height $<n$. By the induction hypothesis, this sequent is valid. To show the conclusion $\Gamma, A\land B \seq C$ is valid, assume that $\mathcal{M},w \models_{g_w} \bigwedge\Gamma \land A\land B$. Then $\mathcal{M},w \models_{g_w} A\land B$, which implies $\mathcal{M},w \models_{g_w} A$ and $\mathcal{M},w \models_{g_w} B$. Hence all formulas of the premise antecedent are true at $w$ under $g_w$, so by validity of the premise sequent we obtain $\mathcal{M},w \models_{g_w} C$. Thus the conclusion is valid.
		
		\smallskip
		\noindent\textbf{Rule $\rAndR$:} The premises are $\Gamma \seq A$ and $\Gamma \seq B$, both valid by the induction hypothesis. Assume $\mathcal{M},w \models_{g_w} \bigwedge\Gamma$. By the first premise, $\mathcal{M},w \models_{g_w} A$; by the second, $\mathcal{M},w \models_{g_w} B$. Hence $\mathcal{M},w \models_{g_w} A\land B$, so the conclusion $\Gamma \seq A\land B$ is valid.
		
		\smallskip
		\smallskip\noindent\textbf{Rule $\rAllR$:} The premise is $\Gamma \seq \varphi[y/x]$ with eigenvariable $y \notin \fv(\Gamma) \cup \fv(\forall x\varphi)$. Suppose $\mathcal{M},w \models_{g_w} \bigwedge\Gamma$. Let $d \in D(w)$ be arbitrary and set $g'_w=g_w[d/y]$. Since $y\notin\fv(\Gamma)$, Coincidence gives $\mathcal{M},w\models_{g'_w}\bigwedge\Gamma$, and the induction hypothesis yields $\mathcal{M},w\models_{g'_w}\varphi[y/x]$. If $y=x$, then $g'_w=g_w[d/x]$ and $\varphi[y/x]=\varphi$, so the required instance follows immediately. If $y\neq x$, Lemma~\ref{lem:sem-subst}(b) gives
		\[
		\mathcal M,w\models_{g'_w[d/x]}\varphi,
		\]
		because $\overline{g'_w}(y)=d$. The eigenvariable condition implies $y\notin\fv(\varphi)$, so $g'_w[d/x]$ and $g_w[d/x]$ agree on $\fv(\varphi)$; Lemma~\ref{lem:sem-subst}(a) therefore gives $\mathcal M,w\models_{g_w[d/x]}\varphi$. Since $d$ was arbitrary, $\mathcal M,w\models_{g_w}\forall x\varphi$.
		
		\smallskip
		\noindent\textbf{Rule $\rAllL$:} The premise is $\Gamma, \forall x\varphi, \varphi[t/x] \seq C$, derived with a smaller height and hence valid by the induction hypothesis. Assume $\mathcal{M},w \models_{g_w} \bigwedge\Gamma \wedge \forall x\varphi$. In particular $\mathcal{M},w \models_{g_w} \forall x\varphi$. Let $d = \overline{g_w}(t)$; regardless of whether or not it is a variable or a constant, by term evaluation clause and constant rigidity, $d \in D(w)$. From the truth condition of $\forall$, $\mathcal{M},w \models_{g_w[d/x]} \varphi$, which is equivalent to $\mathcal{M},w \models_{g_w} \varphi[t/x]$ by Lemma~\ref{lem:sem-subst}(b). Hence all formulas in the antecedent of the premise are true at $w$ under $g_w$, and the induction hypothesis yields $\mathcal{M},w \models_{g_w} C$. So the conclusion $\Gamma, \forall x\varphi \seq C$ is valid.
		
		\smallskip
		\noindent\textbf{Rule $\rDiamK$:} The premise $\varphi \seq \psi$ is valid by the induction hypothesis. Assume $\mathcal{M},w \models_{g_w} \bigwedge\Gamma \wedge \Diamond\varphi$. Then $\mathcal{M},w \models_{g_w} \Diamond\varphi$, so there exists a world $v \in W$ with $w R v$ and $\mathcal{M},v \models_{g_w} \varphi$. Since domains expand ($D(w) \subseteq D(v)$), $g_w$ is also a valid $v$-assignment. The premise $\varphi \seq \psi$ is valid in all worlds under all assignments, so in particular $\mathcal{M},v \models_{g_w} \psi$. Therefore $\mathcal{M},w \models_{g_w} \Diamond\psi$, and the conclusion $\Gamma, \Diamond\varphi \seq \Diamond\psi$ is valid.
		
		\smallskip
		\noindent\textbf{Rule $\rDiamFour$:} The premise $\varphi \seq \Diamond\psi$ is valid by the induction hypothesis. Assume $\mathcal{M},w \models_{g_w}\bigwedge\Gamma \wedge \Diamond\varphi$. Then $\mathcal{M},w \models_{g_w} \Diamond\varphi$, so there exists $v \in W$ with $w R v$ and $\mathcal{M},v \models_{g_w} \varphi$. By the premise, $\mathcal{M},v \models_{g_w} \Diamond\psi$, so there exists $u \in W$ with $v R u$ and $\mathcal{M},u \models_{g_w} \psi$. By transitivity of $R$, $w R u$, hence $\mathcal{M},w \models_{g_w} \Diamond\psi$. The conclusion $\Gamma, \Diamond\varphi \seq \Diamond\psi$ is valid.
		
		\smallskip
		All cut-free rules have been verified. The induction is complete.
	\end{proof}
	
	\subsection{Completeness relative to the Hilbert system\\}
	\label{subsec:complete-hilbert}
	
	We show that every $\QRC$ Hilbert axiom and rule is derivable (using $\rCut$) in $\GQRC$. Combined with Theorem~\ref{thm:cut}, this gives cut-free completeness.
	
	\begin{theorem}[Completeness]
		\label{thm:complete}
		If $\varphi\vdash_{\QRC}\psi$, then $\vdash \varphi\seq\psi$ in $\GQRC$.
	\end{theorem}
	
	\begin{proof}
		We show that each Hilbert axiom and rule of $\QRC$ (Definition~\ref{def:QRC}) is derivable in $\GQRC$. Whenever a cut-free structural lemma is used below, Cut Admissibility is first applied to the given premise derivation.
		
		\medskip
		\noindent\textbf{(i-a)} $\varphi\Hderiv\top$: derivation $\varphi\seq\top$
		by $\rTopR$.
		
		\noindent\textbf{(i-b)} $\varphi\Hderiv\varphi$: derivation $\varphi\seq\varphi$
		by Lemma~\ref{lem:gen-id} (General Identity).
		
		\noindent\textbf{(ii-a)} $\varphi\wedge\psi\Hderiv\varphi$:
		\[
		\infer[\rAndL]{\varphi\wedge\psi\seq\varphi}
		{\infer[\text{Wk. Lemma~\ref{lem:weak}}]{\varphi,\psi\seq\varphi}
			{\infer[\text{Gen.\ Id.}]{\varphi\seq\varphi}{}}}
		\]
		
		\noindent\textbf{(ii-b)} $\varphi\wedge\psi\Hderiv\psi$: symmetric.
		
		\noindent\textbf{(iii)} If $\varphi\seq\psi$ and $\varphi\seq\gamma$, then $\varphi\seq\psi\wedge\gamma$ by $\rAndR$.
		
		\noindent\textbf{(iv) Transitivity} if $\varphi\seq\psi$ and $\psi\seq\gamma$ are derivable, then $\varphi\seq\gamma$ is their $\rCut$.
		
		\noindent\textbf{(v) Monotonicity:} if $\varphi\seq\psi$, then $\Diamond\varphi\seq\Diamond\psi$ by $\rDiamK$ with empty context.
		
		\noindent\textbf{(vi)} $\Diamond\Diamond\varphi\Hderiv\Diamond\varphi$:
		\[
		\infer[\rDiamFour]{\Diamond\Diamond\varphi\seq\Diamond\varphi}
		{\infer[\text{Gen.\ Id.}]{\Diamond\varphi\seq\Diamond\varphi}{}}
		\]
		
		\noindent\textbf{(vii) Quantifier--modality axiom:}
		$\Diamond\forall x\varphi\Hderiv\forall x\Diamond\varphi$. Choose $y\notin\mathrm{Var}(\varphi)\cup\{x\}$ fresh:
		\[
		\infer[\rAllR]{{\Diamond\forall x\varphi\seq\forall x\Diamond\varphi}}
		{\infer[\rDiamK]{\Diamond\forall x\varphi\seq\Diamond\varphi[y/x]}
			{\infer[\rAllL]{\forall x\varphi\seq\varphi[y/x]}
				{\infer[\text{Gen.\ Id.}]{\forall x\varphi,\varphi[y/x]\seq\varphi[y/x]}{}}}}
		\]
		Note $\Diamond\varphi[y/x] = (\Diamond\varphi)[y/x]$. Then by $\rAllR$ with eigenvariable $y$ where $y\notin\fv(\Diamond\forall x\varphi)$: $\Diamond\forall x\varphi\seq\forall x\Diamond\varphi$.
		
		\noindent\textbf{(viii) Generalization:} Suppose $\varphi\seq\psi$ is derivable and $x\notin\fv(\varphi)$. Since $\psi[x/x]=\psi$ and $x\notin\fv(\varphi)\cup\fv(\forall x\psi)$, an application of $\rAllR$ with eigenvariable $x$ yields $\varphi\seq\forall x\psi$.
		
		\noindent\textbf{(ix) Instantiation left:} if $\varphi[t/x]\seq\psi$, then $\forall x\varphi\seq\psi$:
		\[
		\infer[\rAllL]{\forall x\varphi\seq\psi}
		{{\forall x\varphi,\varphi[t/x]\seq\psi}}
		\]
		where $\forall x\varphi, \varphi[t/x]\seq\psi$ is obtained from $\varphi[t/x]\seq\psi$ by Lemma~\ref{lem:weak}.
		
		\noindent\textbf{(x) Substitution:} if $\varphi\seq\psi$, then $\varphi[t/x]\seq\psi[t/x]$ by Lemma~\ref{lem:subst}.
		
		\noindent\textbf{(xi) Eigenvariable rule:}
		Suppose $\varphi[c/x] \seq \psi[c/x]$ is derivable in $\GQRC$, where $c$ is a constant that does \emph{not} occur in $\varphi$ or $\psi$. By Theorem~\ref{thm:cut}, take a cut-free derivation and apply Lemma~\ref{lem:const-subst} (Constant Substitution), replacing the constant $c$ with the variable $x$. This yields a derivation of $(\varphi[c/x])[x/c] \seq (\psi[c/x])[x/c]$ of the same height. Because $c$ does not occur in the original formulas $\varphi$ and $\psi$, a straightforward structural induction on formulas, using the capture-avoiding clauses in Definition~\ref{subsdef}, gives
		\[
		(\varphi[c/x])[x/c] \equiv_\alpha \varphi,
		\qquad
		(\psi[c/x])[x/c] \equiv_\alpha \psi.
		\]
		Because our calculus identifies formulas and sequents up to $\alpha$-equivalence, the derived sequent is $\varphi \seq \psi$.
	\end{proof}
	
	Having shown that the sequent rules capture the Hilbert axioms, we now check their admissibility:
	
	\begin{lemma}[Hilbert admissibility of sequent rules]
		\label{lem:hilbert-admiss}
		If $\vdash \Gamma \seq C$ in $\GQRC$, then $\bigwedge\Gamma \Hderiv C$ in $\QRC$ (where $\bigwedge\emptyset := \top$).
	\end{lemma}
	
	\begin{proof}
		By Theorem~\ref{thm:cut} (Cut Admissibility), every derivable sequent possesses a cut-free derivation. Hence it suffices to prove the statement for cut-free derivations. We proceed by induction on the height $n$ of a cut-free derivation. The induction hypothesis is that for every sequent $\Delta \seq D$ that has a cut-free derivation of height strictly less than $n$, we have $\bigwedge\Delta \Hderiv D$ in $\QRC$. We now examine the last rule of the derivation.
		
		\smallskip
		\noindent\textbf{Axiom $\rId$:} The sequent is $\Gamma', P \seq P$ with
		$P$ atomic. Then $\bigwedge\Gamma' \wedge P \Hderiv P$ is an instance of Hilbert axiom (ii-b).
		
		\smallskip
		\noindent\textbf{Axiom $\rTopR$:} The sequent is $\Gamma \seq \top$. Then $\bigwedge\Gamma \Hderiv \top$ is an instance of Hilbert axiom (i-a). 
		
		\smallskip
		\noindent\textbf{Rule $\rAndL$:} The premise is $\Gamma', A, B \seq C$, derived with height $< n$. By the induction hypothesis, $\bigwedge\Gamma' \wedge A \wedge B \Hderiv C$. From Hilbert axioms (ii-a) and (ii-b) and rule (iii) we obtain $\bigwedge\Gamma' \wedge (A \land B) \Hderiv \bigwedge\Gamma' \wedge A \wedge B$. Transitivity (iv) then yields $\bigwedge\Gamma' \wedge (A \land B) \Hderiv C$.
		
		\smallskip
		\noindent\textbf{Rule $\rAndR$:} The premises are $\Gamma \seq A$ and $\Gamma \seq B$, both derived with height $< n$. By the induction hypothesis, $\bigwedge\Gamma \Hderiv A$ and $\bigwedge\Gamma \Hderiv B$. Rule (iii) gives $\bigwedge\Gamma \Hderiv A \land B$, i.e.\ $\bigwedge\Gamma \Hderiv (A \land B)$.
		
		\smallskip
		\noindent\textbf{Rule $\rAllL$:} The premise is $\Gamma', \forall x\varphi, \varphi[t/x] \seq C$, derived with height $< n$. By the induction hypothesis, $\bigwedge\Gamma' \land \forall x\varphi \land \varphi[t/x] \Hderiv C$. Choose an $\alpha$-equivalent representative $\forall z\varphi'$ of $\forall x\varphi$ such that $z$ and every bound variable of $\varphi'$ avoid $\fv(t)$. Then $t$ is free for $z$ in $\varphi'$, and capture-avoiding substitution gives $\varphi'[t/z]\equiv_\alpha\varphi[t/x]$. Hilbert rule (ix), read under the standing $\alpha$-conversion convention, therefore yields $\forall x\varphi \Hderiv \varphi[t/x]$. Using (ii-a), (ii-b), and (iii), we obtain $\bigwedge\Gamma' \land \forall x\varphi \Hderiv \bigwedge\Gamma' \land \forall x\varphi \land \varphi[t/x]$. Transitivity (iv) with the induction hypothesis then yields $\bigwedge\Gamma' \land \forall x\varphi \Hderiv C$.
		
		\smallskip
		\noindent\textbf{Rule $\rAllR$:} The premise is $\Gamma \seq \varphi[y/x]$ with eigenvariable $y \notin \fv(\Gamma,\forall x\varphi)$. The induction hypothesis gives $\bigwedge\Gamma \Hderiv \varphi[y/x]$. Since $y \notin \fv(\bigwedge\Gamma)$, Hilbert rule (viii) yields
		\[
		\bigwedge\Gamma \Hderiv \forall y(\varphi[y/x]).
		\]
		The latter formula is $\alpha$-equivalent to $\forall x\varphi$, so the standing convention gives $\bigwedge\Gamma \Hderiv \forall x\varphi$.
		
		\smallskip
		\noindent\textbf{Rule $\rDiamK$:}
		The premise is $\varphi \seq \psi$, derived with height $< n$. By the induction hypothesis, $\varphi \Hderiv \psi$. Hilbert rule (v) yields $\Diamond\varphi \Hderiv \Diamond\psi$. From axiom (ii-b) we have $\bigwedge\Gamma \land \Diamond\varphi \Hderiv \Diamond\varphi$. Applying rule (iv) to these two statements gives $\bigwedge\Gamma \land \Diamond\varphi \Hderiv \Diamond\psi$.

		\smallskip
		\noindent\textbf{Rule $\rDiamFour$:}
		The premise is $\varphi \seq \Diamond\psi$, derived with height $< n$. By the induction hypothesis, $\varphi \Hderiv \Diamond\psi$. Rule (v) gives $\Diamond\varphi \Hderiv \Diamond\Diamond\psi$; axiom (vi) gives $\Diamond\Diamond\psi \Hderiv \Diamond\psi$. Transitivity (iv) yields $\Diamond\varphi \Hderiv \Diamond\psi$. As in the $\rDiamK$ case, we use (ii-b): $\bigwedge\Gamma \land \Diamond\varphi \Hderiv \Diamond\varphi$, and transitivity gives $\bigwedge\Gamma \land \Diamond\varphi \Hderiv \Diamond\psi$, as required.
		
		\smallskip
		All cut-free rules have been treated, so the induction is complete.
	\end{proof}
	
	\begin{corollary}[Equivalence]
		\label{cor:equiv}
		For any formulas $\varphi,\psi$,
		\[
		\varphi\vdash_{\QRC}\psi \text{ in } \QRC
		\quad\Longleftrightarrow\quad
		\vdash \varphi\seq\psi \text{ in } \GQRC .
		\]
	\end{corollary}
	
	\begin{proof}
		$(\Rightarrow)$: Theorem~\ref{thm:complete} directly shows that every Hilbert derivation of $\varphi\vdash_{\QRC}\psi$ can be simulated in $\GQRC$. Hence $\vdash \varphi\seq\psi$.
		
		$(\Leftarrow)$: Assume $\vdash \varphi\seq\psi$ in $\GQRC$ (possibly using cuts). By Theorem~\ref{thm:cut} (Cut Admissibility), this sequent also has a cut-free derivation. Applying Lemma~\ref{lem:hilbert-admiss} to that cut-free derivation yields $\varphi\Hderiv\psi$ in $\QRC$.
	\end{proof}

	\subsection{Finite Completeness and the Finite Model Property\\}
	\label{subsec:finite-completeness}
	
	As we had simplified the semantics of Borges and Joosten \cite{BorgesJoosten2020}, we now give a direct completeness proof that simultaneously yields the finite model property. Fix a sequent $\Gamma_0 \seq C_0$. Let $\mathcal{T}_0$ be the finite set of all terms occurring in $\Gamma_0 \cup \{C_0\}$. Let $Q$ be the total number of occurrences of $\forall$ in $\Gamma_0 \cup \{C_0\}$. We work with a countably infinite reservoir of individual variables from which all fresh variables needed in the construction are drawn. Freshness is global: every newly chosen variable is absent from all formulas and domains constructed earlier.
	
	\textit{Left saturation} expands the antecedent of a sequent by adding every formula that the left logical rules could extract from it.
	
	\begin{definition}[Left saturation]
		\label{def:left-sat}
		Let $\Theta$ be a finite multiset of formulas and $T$ a finite set of terms. The \emph{left saturation} $\Theta^{\dagger}$ (over $T$) is the least set of $\alpha$-equivalence classes of formulas that contains the support of $\Theta$ and is closed under:
		
		\begin{itemize}
			\item if $A \land B \in \Theta^{\dagger}$, then $A, B \in \Theta^{\dagger}$;
			\item if $\forall x A \in \Theta^{\dagger}$ and $t \in T$, then $A[t/x] \in \Theta^{\dagger}$.
		\end{itemize}
		When used as an antecedent, $\Theta^{\dagger}$ is represented by one formula from each class, with multiplicity one. The original formulas are retained up to $\alpha$-equivalence. This set is finite. Indeed, every generated formula is, up to $\alpha$-equivalence, of the form $B\sigma$, where $B$ is a subformula of a formula in $\Theta$ and $\sigma$ is a partial substitution whose domain is contained in $\fv(B)$ and whose values lie in $T$. The conjunction clause passes to a conjunct, and the quantifier clause extends such a partial substitution by sending the newly exposed variable to $t\in T$. Since there are only finitely many subformulas $B$, and at most $(|T|+1)^{|\fv(B)|}$ such partial substitutions for each $B$, there are at most
		\[
		\sum_B (|T|+1)^{|\fv(B)|}
		\]
		distinct classes. Hence $\Theta^{\dagger}$ is finite, including when $T=\varnothing$.
	\end{definition}
	
	\begin{lemma}[Left-Saturation Entailment]
		\label{lem:left-sat-equiv}
		Let $\Theta$ be a finite multiset of formulas, $T$ a finite set of terms, and $\Theta^{\dagger}$ the left saturation of $\Theta$ over $T$ as in Definition~\ref{def:left-sat}. If $\vdash \Theta^{\dagger} \seq F$ then $\vdash \Theta \seq F$.
	\end{lemma}
	\begin{proof}
		By Cut Admissibility (Theorem~\ref{thm:cut}), replace the assumed derivation of $\Theta^{\dagger}\seq F$ by a cut-free derivation. Every transformation below therefore remains within the cut-free calculus, and Theorem~\ref{thm:contr} is applicable. Because $\Theta^{\dagger}$ is finite, choose a finite staged construction from the support of $\Theta$ by operations that add $A,B$ alongside an existing $A\land B$, or add $\varphi[t/x]$ alongside an existing $\forall x\varphi$. Remove these additions in reverse order. The carrier formula responsible for each addition is then still present, and the following local reductions apply.
		
		Suppose a reverse stage removes one or both conjuncts generated
		from $A\land B\in\Delta$. Use height-preserving weakening to supply
		any missing occurrence of $A$ or $B$, obtaining
		$\Delta,A,B\seq F$. Writing $\Delta=\Delta_0,A\land B$, an application
		of $\rAndL$ followed by contraction yields $\Delta\seq F$.
		
		Suppose a derivation of $\Delta,\varphi[t/x]\seq F$ is given, where $\forall x\varphi\in\Delta$. Write $\Delta=\Delta_0,\forall x\varphi$. Then the given end-sequent is $\Delta_0,\forall x\varphi,\varphi[t/x]\seq F$. Applying $\rAllL$ yields $\Delta_0,\forall x\varphi\seq F$, which is exactly $\Delta\seq F$.
		
		Applying these reductions finitely many times gives a derivation from the support of $\Theta$. Height-preserving weakening then restores any additional multiplicities present in the original multiset $\Theta$, proving $\Theta\seq F$.
	\end{proof}
	For a finite multiset $\Theta$, write
	\[
	\mathrm{md}(\Theta)=\max\{\mathrm{md}(A):A\in\Theta\},
	\qquad
	\mathrm{md}(\varnothing)=0.
	\]
	
	\begin{definition}[Canonical model]
		\label{def:canonical}
		We build a finite Kripke model $\mathcal{M} = (W, R, D, I)$ inductively, starting from the sequent $\Gamma_0 \seq C_0$. Since the set of syntactic variables and constants is countable, fix an injection of all terms into the ambient universe $\mathcal U$ and, to simplify notation, identify each term with its image. Thus every set of terms used below is literally a set of individuals.
		
		Every world $w\in W$ will be equipped with a finite multiset $\Theta_w$, called the \emph{root antecedent}, and a finite local domain $D(w)$. From these we form the \emph{saturated antecedent} $\Gamma_w^*=(\Theta_w)^{\dagger}$ over $D(w)$ (Definition~\ref{def:left-sat}).
		
		Set $\Theta_{w_0}=\Gamma_0$ and choose a set $Z_{w_0}$ of $Q+1$ variables fresh from $\Gamma_0\cup\{C_0\}$. Define
		\[
		D(w_0)=\mathcal T_0\cup Z_{w_0},
		\]
		and saturate the antecedent to obtain $\Gamma_{w_0}^*$. The local domain is never empty.
		
		\smallskip\noindent\textbf{Successor worlds}: For each formula $\Diamond A\in\Gamma_w^*$, create a new world $v$ with root antecedent $\Theta_v=\{A\}$. Choose a set $Z_v$ of $Q+1$ fresh variables from the global reservoir, and define
		\[
		D(v)=D(w)\cup Z_v.
		\]
		Saturate $\Theta_v$ over $D(v)$ to obtain $\Gamma_v^*$.
		
		Whenever a successor $v$ is built from a parent $w$, we add the pair $(w,v)$ to $R_0$. After the inductive construction finishes, let $W$ be the set of all worlds obtained. The accessibility relation $R$ of the model is the strict transitive closure $R_0^+$.
		
		Every successor world $v$ generated from a world $w$ is initialized with a root antecedent $\Theta_v = \{A\}$, triggered by some diamond formula $\Diamond A$ present in the saturated multiset $\Gamma_w^*$. We can track the maximum modal depth of these multisets to guarantee termination. Recall that $\Gamma_w^*$ is constructed from the root $\Theta_w$ using only non-modal saturation rules. Because non-modal rules never increase modal depth, the maximum modal depth in $\Gamma_w^*$ is bounded by the maximum modal depth in $\Theta_w$. Since $\Diamond A \in \Gamma_w^*$, its depth obeys this bound. When we strip the modality to form the new root, its depth strictly drops: $\mathrm{md}(A) < \mathrm{md}(\Diamond A)$. Chaining these facts together, the modal depth of the new root $\Theta_v$ is strictly smaller than the modal depth of the parent root $\Theta_w$. 
		
		Consequently, the maximum modal depth of the root antecedents strictly decreases along every step of the generating relation $R_0$. Since modal depth is a non-negative integer, any chain of successors starting from the initial world $w_0$ must terminate in a finite number of steps, bounded by the maximum modal depth of the original end-sequent $\Gamma_0 \seq C_0$. Furthermore, because the saturated multiset $\Gamma_w^*$ is finite, each world can trigger only finitely many successors. Since the generated frame is a tree of finite depth with finite branching, the total set of generated worlds $W$ is strictly finite. Finally, because the accessibility relation $R$ is defined as the transitive closure of $R_0$, and domains expand along every direct successor step ($D(w) \subseteq D(v)$ by construction), transitivity natively guarantees that the expanding domains condition holds globally.
		
		To ensure truth in our model perfectly mirrors the syntactical structure of our derivations, we construct a \emph{term model}. An atomic formula is true under the canonical assignment precisely when the atom obtained by evaluating its arguments occurs in the corresponding saturated antecedent.
		
		\smallskip\noindent\textbf{Interpretation of constants.}
		Let $\Sigma$ be the full signature. For every constant $c \in \Sigma$ that belongs to the active term set $\mathcal{T}_0$, we set $I(c) = c$, so it denotes itself. For all other constants $c \in \Sigma$, we choose an arbitrary but fixed element of $D(w_0)$ and set $I(c)$ to that element. Because this assigned element belongs to the root domain $D(w_0)$, and domains monotonically expand along $R$, it belongs to every $D(w)$. Hence, the constant rigidity condition $I(c) \in D(w)$ is globally satisfied for all worlds $w \in W$.
		
		No new constant is introduced during the construction: every successor root is the scope of a diamond already present in a parent saturation, and left saturation instantiates only with terms from the current domain, whose constant members all originate in $\mathcal T_0$. Consequently, every constant occurring in any $\Gamma_w^*$ belongs to $\mathcal T_0$ and therefore denotes itself.
		
		\smallskip\noindent\textbf{Valuation and assignments.} For each $n$-ary relation symbol $S$ and each world $w$, the interpretation $I_w(S)$ is defined strictly syntactically:
		\[
		I_w(S) = \{\, (t_1,\dots,t_n) \in D(w)^n \mid S(t_1,\dots,t_n) \in \Gamma_w^* \,\}.
		\]
		To evaluate variables, we fix an arbitrary default element $d_w \in D(w)$ for each world $w$. The $w$-assignment $g_w : \mathsf{Var} \to D(w)$ is then defined as the identity on the local domain:
		\[
		g_w(x) =
		\begin{cases}
			x, & \text{if } x \in D(w), \\
			d_w, & \text{otherwise}.
		\end{cases}
		\]
		By construction, all variables that occur free in the formulas of $\Gamma_w^*$ are contained in $D(w)$. Consequently, they evaluate exactly to themselves under $g_w$. This completes the definition of the finite $\QRC$-model $\mathcal{M}$.
	\end{definition}
	
	In the remainder of the canonical-model argument, $w\models A$ abbreviates $\mathcal M,w\models_{g_w}A$.
	
	\begin{lemma}[Antecedent Truth]
		\label{lem:ant-truth}
		For every world $w\in W$ and every formula $A\in\Gamma_w^*$, $\mathcal M,w\models_{g_w}A$.
	\end{lemma}
	\begin{proof}
		By induction on the logical complexity $|A|$. The induction hypothesis is that the claim holds for all worlds and all formulas of strictly smaller complexity.
		
		\begin{itemize}
			\item \emph{$A=\top$.} This is immediate from the satisfaction clause for $\top$.
			\item \emph{Atomic $S(t_1,\dots,t_n)$.}
			Assume $S(t_1,\dots,t_n) \in \Gamma_w^*$. For every term $t_i$ we have $\overline{g_w}(t_i) = t_i$: variables occurring in $\Gamma_w^*$ are in $D(w)$ and $g_w$ maps them to themselves; constants occurring in $\Gamma_w^*$ belong to $\mathcal T_0$ and are interpreted as themselves. Hence $(\overline{g_w}(t_1),\dots,\overline{g_w}(t_n)) = (t_1,\dots,t_n)$. By definition of $I_w(S)$ and the assumption, the tuple $(t_1,\dots,t_n)$ belongs to $I_w(S)$. Therefore the evaluated tuple is in $I_w(S)$, which means $\mathcal{M},w \models S(t_1,\dots,t_n)$.
			\item \emph{$A = B \land C$.} 
			Left saturation guarantees $B, C \in \Gamma_w^*$. Applying the induction hypothesis to $B$ and $C$ (each of smaller complexity) yields $w \models B$ and $w \models C$, so $w \models B \land C$.
			\item \emph{$A = \forall x B$.} 
			For every $d \in D(w)$, left saturation placed $B[d/x] \in \Gamma_w^*$. Since $|B[d/x]| < |\forall x B|$, the induction hypothesis yields $\mathcal{M}, w \models_{g_w} B[d/x]$. Every member $d$ of $D(w)$ is a syntactic term and denotes itself at $w$: variables by the definition of $g_w$, and constants by the definition of $I$. Lemma~\ref{lem:sem-subst}(b) therefore gives $\mathcal{M}, w \models_{g_w} B[d/x]$ iff $\mathcal{M}, w \models_{g_w[d/x]} B$. Thus $w \models_{g_w} \forall x B$.
			\item \emph{$A = \Diamond B$.} By construction of the canonical model, the presence of $\Diamond B$ in $\Gamma_w^*$ guarantees the existence of a world $v$ such that $w R_0 v$ and $\Theta_v = \{B\}$. Consequently, $B \in \Gamma_v^*$, as left saturation retains the root formulas. Because $|B| < |\Diamond B|$ and $B \in \Gamma_v^*$, the induction hypothesis strictly applies to formula $B$ at world $v$, yielding $v \models_{g_v} B$. We must show $v \models_{g_w} B$, as the truth condition of $\Diamond$ at $w$ uses the assignment $g_w$. All free variables of $B$ occur in $\Diamond B$ and therefore lie in $D(w)$ as $\Diamond B \in \Gamma_w^*$. Since $D(w) \subseteq D(v)$, both $g_w$ and $g_v$ act as the identity on these variables, so they agree on the free variables of $B$. Consequently $v \models_{g_v} B$ implies $v \models_{g_w} B$. Together with $w R v$ (as $R_0 \subseteq R$) this gives $w \models_{g_w} \Diamond B$.
		\end{itemize}
	\end{proof}
	
	For a world $w$, call $F$ a \emph{$D(w)$-instance} of a formula $G$ if, up to $\alpha$-equivalence, $F=G\sigma$ for a finite simultaneous capture-avoiding substitution $\sigma$ such that $\sigma(x)\in D(w)$ for every $x\in\fv(G)$, including when $\sigma(x)=x$.
	
	\begin{lemma}[Truth Implies Provability]
		\label{lem:truth-prov}
		Let $w\in W$ and let $F$ be a $D(w)$-instance of a subformula of $C_0$. Suppose that $k_w(F)+n(F)\le Q$, where $k_w(F)=|Z_w\cap\fv(F)|$ and $n(F)$ is the number of $\forall$-symbols in $F$. If $\mathcal{M},w\models_{g_w}F$ then $\vdash\Gamma_w^*\seq F$.
	\end{lemma}
	
	\begin{proof}
		We proceed by an outer induction on $\mathrm{md}(\Theta_w)$ and an inner induction on $|F|$. The outer induction hypothesis is that for every world $u$ with $\mathrm{md}(\Theta_u)<\mathrm{md}(\Theta_w)$, the statement of the lemma holds at $u$ for all formulas satisfying its conditions. The inner induction hypothesis is that for the current world $w$, the statement holds for all $G$ with $|G|<|F|$ that satisfy $k_w(G)+n(G)\le Q$. As established in the canonical model construction, $wR_0u$ implies $\mathrm{md}(\Theta_u)<\mathrm{md}(\Theta_w)$.
		
		\medskip\noindent\textbf{Preservation of the hypotheses.} The numerical hypothesis gives $k_w(F)+n(F)\le Q$. In the conjunction case, the inner induction passes to a conjunct $G$ with $\fv(G)\subseteq\fv(F)$ and $n(G)\le n(F)$, so $k_w(G)+n(G)\le Q$. In the universal case, it passes to $A[z/x]$ with $z\in Z_w$ fresh for $\forall xA$; this decreases $n$ by exactly $1$ and increases $k_w$ by at most $1$. In the modal case, the outer induction is applied at a direct successor, where the numerical condition is checked explicitly below.
		
		The syntactic hypothesis is preserved as well. If a conjunction or a diamond is a $D(w)$-instance of a subformula of $C_0$, then each immediate component to which the induction passes is a $D(w)$-instance of the corresponding immediate subformula; no binder is crossed in either case. If $\forall xA$ is a $D(w)$-instance and $z\in D(w)$, then $A[z/x]$ is a $D(w)$-instance of the corresponding matrix subformula, obtained after standardizing binders apart by extending the witnessing substitution with $x\mapsto z$. Finally, whenever $D(w)\subseteq D(u)$, every $D(w)$-instance is also a $D(u)$-instance. Hence every invocation of either induction hypothesis satisfies both the syntactic and the numerical assumptions of the lemma.
		
		\medskip
		\noindent\emph{Case $F=\top$.} The sequent $\Gamma_w^*\seq\top$ is an instance of $\rTopR$.
		
		\smallskip\noindent\emph{Case $F$ atomic.}
		Write $F=S(t_1,\dots,t_n)$. From $w\models_{g_w}F$ we obtain $(\overline{g_w}(t_1),\dots,\overline{g_w}(t_n))\in I_w(S)$. Every variable among the $t_i$ belongs to $D(w)$ by the hypothesis on $F$ and hence evaluates to itself. Every constant among the $t_i$ belongs to $D(w)$; by the construction all constant members of any local domain originate in $\mathcal T_0$, and so are interpreted as themselves. Hence the tuple $(t_1,\dots,t_n)$ belongs to $I_w(S)$, which by definition means $F\in\Gamma_w^*$. The axiom $\rId$ then gives $\vdash\Gamma_w^*\seq F$.
		
		\smallskip\noindent\emph{Case $F=A\land B$.}
		$w\models_{g_w}A\land B$ implies $w\models_{g_w}A$ and $w\models_{g_w}B$. Since $|A|,|B|<|A\land B|$, and the invariants $k_w(A)+n(A)\le Q$ and $k_w(B)+n(B)\le Q$ follow from the preservation clause, the inner induction hypothesis yields $\vdash\Gamma_w^*\seq A$ and $\vdash\Gamma_w^*\seq B$. Applying $\rAndR$ gives $\vdash\Gamma_w^*\seq A\land B$.
		
		\smallskip\noindent\emph{Case $F=\forall xA$.} Assume $w\models_{g_w}\forall xA$; then $w\models_{g_w}A[d/x]$ for every $d\in D(w)$. From the hypothesis $k_w(F)+n(F)\le Q$ and $n(F)\ge 1$ we obtain $k_w(F)\le Q-1$. Because $k_w(F) = |Z_w \cap \fv(F)|$ specifically counts how many variables from the reservoir $Z_w$ appear free in $F = \forall xA$, our established bound $k_w(F) \le Q-1$ guarantees that at most $Q-1$ of the $Q+1$ available variables in $Z_w$ are used, forcing the set difference $Z_w \setminus \fv(\forall xA)$ to be nonempty; choose $z$ from it. Then $z\notin\fv(\forall xA)$ by choice, and $z\notin\fv(\Theta_w)$ because $Z_w\cap\fv(\Theta_w)=\emptyset$ by construction.
		
		Since $z\in D(w)$, the truth of $\forall xA$ yields $w\models_{g_w}A[z/x]$. The formula $A[z/x]$ satisfies $k_w(A[z/x])\le k_w(F)+1$ and $n(A[z/x])=n(F)-1$, so $k_w(A[z/x])+n(A[z/x])\le Q$. Moreover $|A[z/x]|<|\forall xA|$. The inner induction hypothesis therefore gives $\vdash\Gamma_w^*\seq A[z/x]$. Lemma~\ref{lem:left-sat-equiv} shrinks the antecedent to the root: $\vdash\Theta_w\seq A[z/x]$. The eigenvariable condition $z\notin\fv(\Theta_w)\cup\fv(\forall xA)$ holds, so $\rAllR$ yields $\vdash\Theta_w\seq\forall xA$. Repeated contraction reduces the antecedent to the support of $\Theta_w$, which is contained in $\Gamma_w^*$, and Weakening then yields $\vdash\Gamma_w^*\seq\forall xA$.
		
		\smallskip\noindent\emph{Case $F=\Diamond A$.} Assume $w\models_{g_w}\Diamond A$. By the truth condition there exists a world $v$ with $w R v$ and $v\models_{g_w}A$. Because $\fv(A)\subseteq\fv(\Diamond A)\subseteq D(w)$ and domains expand, $g_w$ and $g_v$ coincide on $\fv(A)$, giving $v\models_{g_v}A$. Since $R=R_0^+$, there is a direct $R_0$-successor $u$ of $w$ such that either $u=v$ or $u R v$.
		
		The successor world $u$ was spawned specifically to satisfy some diamond formula $\Diamond B \in \Gamma_w^*$, meaning its root antecedent was initialized as $\Theta_u = \{B\}$. Since $\Diamond A$ is a $D(w)$-instance of a subformula of $C_0$, its scope $A$ is a $D(w)$-instance of the corresponding scope subformula; and because $D(w)\subseteq D(u)$, both $A$ and $\Diamond A$ satisfy the required local-instance condition at $u$. We must also establish the numerical invariant $k_u(A)=0$. Because the formula $\Diamond A$ is the subject of the induction at world $w$, its free variables are strictly bounded by the parent domain: $\fv(A) \subseteq \fv(\Diamond A) \subseteq D(w)$. When constructing the successor world $u$, the fresh reservoir $Z_u$ is drawn to be explicitly disjoint from all previously constructed sets, including the parent domain $D(w)$. Because all free variables of $A$ are trapped within $D(w)$, the disjointness $Z_u \cap D(w) = \emptyset$ forces $Z_u \cap \fv(A) = \emptyset$ and $Z_u \cap \fv(\Diamond A) = \emptyset$. Thus, $k_u(A) = 0$ and $k_u(\Diamond A) = 0$.

		\begin{itemize}
			\item \emph{If $u=v$:} We have $u\models_{g_u}A$. Since $k_u(A)=0$ and $n(A)\le Q$, the invariant holds. Because $w R_0 u$, the outer induction hypothesis triggers for formula $A$, yielding $\vdash\Gamma_u^*\seq A$. Lemma~\ref{lem:left-sat-equiv} shrinks the antecedent to the root, giving $\vdash\Theta_u\seq A$. Since $\Theta_u = \{B\}$, this is exactly the sequent $\vdash B\seq A$.
			
			\item \emph{If $u R v$:} Since $\fv(A)\subseteq D(w)\subseteq D(u)\subseteq D(v)$, the assignments $g_u$ and $g_v$ agree on $\fv(A)$, so $v\models_{g_u}A$. By the semantic truth condition, this implies $u\models_{g_u}\Diamond A$. As $k_u(\Diamond A)=0$ and $n(\Diamond A)\le Q$, the invariant holds. The outer induction hypothesis applied to formula $\Diamond A$ at world $u$ yields $\vdash\Gamma_u^*\seq\Diamond A$. Lemma~\ref{lem:left-sat-equiv} shrinks this to $\vdash\Theta_u\seq\Diamond A$, which is exactly the sequent $\vdash B\seq\Diamond A$.
		\end{itemize}
		In both cases, we successfully derived either $\vdash B\seq A$ or $\vdash B\seq\Diamond A$. Applying the appropriate modal rule ($\rDiamK$ or $\rDiamFour$) yields $\vdash\Diamond B\seq\Diamond A$. Finally, because $\Diamond B\in\Gamma_w^*$ by definition of the successor $u$, weakening gives the required conclusion $\vdash\Gamma_w^*\seq\Diamond A$.
	\end{proof}
	
	\begin{theorem}[Finite Completeness / Finite Model Property]
		\label{thm:finite-completeness}
		If $\Gamma_0 \seq C_0$ is unprovable in $\GQRC$, then there exists a finite $\QRC$-model that falsifies it.
	\end{theorem}
	\begin{proof}
		Construct the finite model $\mathcal M$ as in Definition~\ref{def:canonical}. By Lemma~\ref{lem:ant-truth}, every formula in $\Gamma_0$ is true at $w_0$ under $g_{w_0}$. Since $\fv(C_0)\subseteq\mathcal T_0\subseteq D(w_0)$, the formula $C_0$ is a $D(w_0)$-instance of itself under the identity substitution. Moreover, $k_{w_0}(C_0)=0$ by freshness of $Z_{w_0}$, and $n(C_0)\leq Q$. If $w_0\models_{g_{w_0}}C_0$, Lemma~\ref{lem:truth-prov} yields $\vdash\Gamma_{w_0}^*\seq C_0$. By Lemma~\ref{lem:left-sat-equiv}, we obtain $\vdash\Gamma_0\seq C_0$, contradicting the assumption. Thus $\mathcal M$ is a finite countermodel to $\Gamma_0\seq C_0$.
	\end{proof}
	
	\begin{corollary}[Decidability]
		For every effectively presented countable signature, $\QRC$ is decidable: there is an algorithm that, given any sequent $\Gamma\seq C$, decides whether it is derivable in $\GQRC$ (equivalently, whether $\bigwedge\Gamma\Hderiv C$ in $\QRC$).
	\end{corollary}
	\begin{proof}
		Here an effective presentation provides decidable symbol equality and arities and an effective finite-string encoding of formulas. The parenthetical equivalence follows from Lemma~\ref{lem:hilbert-admiss} in one direction and, in the other, from Theorem~\ref{thm:complete} followed by Cut Admissibility (Theorem~\ref{thm:cut}) and repeated $\wedge L$-invertibility (or Lemma~\ref{lem:top} when $\Gamma=\varnothing$). Given $\Gamma\seq C$, construct the finite model of Definition~\ref{def:canonical} with $\Gamma_0=\Gamma$ and $C_0=C$. The construction is effective, since left saturation is finite and $\alpha$-equivalence is decidable. By Lemma~\ref{lem:ant-truth} and Soundness, if $\Gamma\seq C$ is derivable then $w_0\models_{g_{w_0}}C$. Conversely, if $w_0\models_{g_{w_0}}C$, Lemmas~\ref{lem:truth-prov} and~\ref{lem:left-sat-equiv} give $\vdash\Gamma\seq C$, since $C$ satisfies the hypotheses of the former at $w_0$. Satisfaction in the constructed finite model is decidable.
	\end{proof}

	\section{Consequences}
	\label{sec:consequences}
	
	For a set of terms $\mathcal T$, a formula $F$ is a
	\emph{$\mathcal T$-instance} of a formula $G$ if
	\[
	F\equiv_\alpha G\sigma
	\]
	for a finite simultaneous substitution $\sigma$ whose nonidentity values
	belong to $\mathcal T$.
	
	\subsection{Subformula Property\\}
	\begin{theorem}[Subformula Property]
		\label{thm:subfml}
		If $\vdash \Gamma \seq C$ in $\GQRC$, then there exists a cut-free derivation $\pi$ of $\Gamma \seq C$ with the following property. Let $T_0$ be the set of all terms occurring in $\Gamma \cup \{C\}$. Let $T_0^* = T_0$ if $T_0$ is non-empty, and $T_0^* = \{v_0\}$ for a fresh variable $v_0$ if $T_0$ is empty. Let $E$ be the finite set of eigenvariables introduced by the applications of $\rAllR$ in $\pi$, and set $\mathcal{T} = T_0^* \cup E$.
		Then:
		\begin{enumerate}
			\item[(i)] every formula occurring in $\pi$ is obtained from a subformula of some formula in $\Gamma \cup \{C\}$ by renaming bound variables and substituting terms from $\mathcal{T}$ for free variables (i.e., is a $\mathcal{T}$-instance of a subformula);
			\item[(ii)] the modal depth of every formula in $\pi$ does not exceed the maximum modal depth of a formula in $\Gamma \cup \{C\}$;
			\item[(iii)] every term used as the witness of an application of $\rAllL$ belongs to $\mathcal{T}$.
		\end{enumerate}
		A cut-free derivation satisfying clauses~(i)--(iii), with $T_0$, $T_0^*$, $E$, and $\mathcal T$ chosen as above, will be called \emph{normalized}.
	\end{theorem}
	
	\begin{proof}
		We first normalize the witness terms. By cut elimination (Theorem~\ref{thm:cut}) choose a cut-free derivation $\pi_0$ of $\Gamma\seq C$. Choose $s^*\in T_0^*$; when $T_0=\varnothing$, take the stipulated variable $v_0$ fresh from the whole derivation.
		
		Let $V_{\mathrm{bad}}$ be the finite set of variable symbols used as $\rAllL$ witnesses in $\pi_0$ that do not belong to $T_0^*$. Apply Lemma~\ref{lem:eigenrename} with $V=V_{\mathrm{bad}}\cup\fv(s^*)$, obtaining a derivation $\pi_1$, and let $E$ be its set of eigenvariables. Every variable witness of $\pi_1$ outside $T_0^*\cup E$ is then a member of $V_{\mathrm{bad}}$: eigenvariable renaming can change such a witness name only when it is an occurrence of the renamed eigenvariable, in which case its new name lies in $E$.
		
		Now replace throughout $\pi_1$ every constant witness outside $T_0^*$ and every variable witness in $V_{\mathrm{bad}}$ by $s^*$, using the rule-by-rule transformations from Lemmas~\ref{lem:const-subst} and~\ref{lem:subst}. Capture avoidance is handled by the local $\alpha$-conversion built into those lemmas, while $E\cap(V_{\mathrm{bad}}\cup\fv(s^*))=\varnothing$ preserves every eigenvariable condition. Hence no eigenvariable is renamed during these replacements and $E$ remains fixed. Every replaced constant is absent from the end-sequent, and every replaced variable has no free occurrence there; hence the substitutions leave the end-sequent unchanged up to $\alpha$-equivalence. If the quantified variable of an $\rAllL$ occurrence is vacuous, its instance formula is independent of the displayed witness, and that rule occurrence is simply relabelled with $s^*$. After finitely many replacements we obtain a cut-free derivation $\pi$ in which every $\rAllL$ witness belongs to $T_0^*\cup E$. Set $\mathcal T=T_0^*\cup E$. This establishes clause~(iii). Bound-variable names introduced solely by $\alpha$-conversion need not lie in $\mathcal T$, because clause~(i) explicitly permits bound-variable renaming.
		
		We now fix the global pool $\mathcal T=T_0^*\cup E$. Read root-first, only the quantifier rules can introduce new free variables, all belonging to $\mathcal T$; hence every free variable occurring in $\pi$ belongs to $\mathcal T$. We prove property~(i) node-wise by induction on subderivation height: for every node of $\pi$, and any choice of syntactic representatives of its sequent, each formula in the subtree above that node is a $\mathcal T$-instance of a subformula of a formula in that node's sequent.
		
		We distinguish cases according to the last rule applied in $\pi$.
		
		\smallskip\noindent\emph{Case $\rId$:} $\pi$ consists of the single sequent $\Gamma', P \seq P$ with $P$ atomic. Every formula in this sequent is a subformula of itself; no substitution is needed.
		
		\smallskip\noindent\emph{Case $\rTopR$:} $\pi$ consists of $\Gamma \seq \top$. The formula $\top$ is a subformula of itself, and no new terms appear.
		
		\medskip
		\smallskip\noindent\emph{Case $\rAndL$:}
		\[
		\infer[\rAndL]{\Gamma', A\land B \seq C}{\Gamma', A, B \seq C}.
		\]
		Let $\pi_1$ be the sub-derivation of the premise, which has strictly smaller height. By the induction hypothesis, every formula occurring in $\pi_1$ is a $\mathcal{T}$-instance of a subformula of some formula in the premise end-sequent $\Gamma' \cup \{A,B,C\}$. Since $A$ and $B$ are proper subformulas of $A\land B$, any $\mathcal{T}$-instance of a subformula of $A$ or $B$ is established as a $\mathcal{T}$-instance of a subformula of $A\land B$. The formulas $\Gamma'$ and $C$ appear unchanged. Thus, every formula in $\pi_1$ is a $\mathcal{T}$-instance of a subformula of the conclusion end-sequent. 
		
		\smallskip\noindent\emph{Case $\rAndR$:}
		\[
		\infer[\rAndR]{\Gamma \seq A\land B}{\Gamma \seq A & \Gamma \seq B}.
		\]
		The argument is similar to $\rAndL$. By the induction hypothesis, every formula in the left premise derivation is a $\mathcal{T}$-instance of a subformula of $\Gamma \cup \{A\}$. Since $A$ is a proper subformula of $A\land B$, all such formulas are established as $\mathcal{T}$-instances of subformulas of $\Gamma \cup \{A\land B\}$. The exact same reasoning holds for the right premise with $B$. 
		
		\smallskip\noindent\emph{Case $\rAllL$:}
		\[
		\infer[\rAllL]{\Gamma', \forall x\varphi \seq C}
		{\Gamma', \forall x\varphi,\; \varphi[t/x] \seq C}.
		\]
		By the induction hypothesis on the strictly shorter premise derivation, every formula in it is a $\mathcal{T}$-instance of a subformula of $\Gamma' \cup \{\forall x\varphi, \varphi[t/x], C\}$. Let $F\equiv_\alpha\psi'\sigma$ be a formula occurring in the premise derivation, where $\psi'$ is a subformula of $\varphi[t/x]$ and the nonidentity values of $\sigma$ belong to $\mathcal T$. Let $\psi$ be the subformula of $\varphi$ at the corresponding occurrence. The free occurrences of any one variable in $\psi$ are either free in $\varphi$ or bound by the same quantifier above $\psi$. Hence $\psi'\equiv_\alpha\psi\rho$ for a simultaneous substitution $\rho$ sending each free variable of $\psi$ either to itself, to $t$, or to the renamed variable of its enclosing binder. Define $\tau(z)=(z\rho)\sigma$ for $z\in\fv(\psi)$ and $\tau(z)=z$ otherwise. After standardizing binders apart, we obtain $F\equiv_\alpha\psi\tau$. For $z\in\fv(\psi)$, if $\tau(z)$ is a variable, it occurs free in $F$ and therefore belongs to $\mathcal T$; if it is a constant, it is either $t$ or a nonidentity value of $\sigma$, and hence again belongs to $\mathcal T$. Consequently $F$ is a $\mathcal T$-instance of the subformula $\psi$ of $\varphi$, and hence of $\forall x\varphi$.
		
		\smallskip\noindent\emph{Case $\rAllR$:}
		\[
		\infer[\rAllR]{\Gamma \seq \forall x\varphi}
		{\Gamma \seq \varphi[y/x]},
		\]
		where $y$ is a fresh eigenvariable recorded in $E \subseteq \mathcal{T}$. The argument is the same as in the $\rAllL$ case, with $t=y$.

		\smallskip\noindent\emph{Case $\rDiamK$:}
		\[
		\infer[\rDiamK]{\Gamma,\Diamond\varphi \seq \Diamond\psi}
		{\varphi \seq \psi}.
		\]
		The premise has smaller height. By the induction hypothesis, every formula in its derivation is a $\mathcal{T}$-instance of a subformula of the premise end-sequent $\{\varphi, \psi\}$. Since $\varphi$ and $\psi$ are proper subformulas of $\Diamond\varphi$ and $\Diamond\psi$, both of which belong to the conclusion, every formula in the upper derivation satisfies the subformula condition with respect to the conclusion. 
		
		\smallskip\noindent\emph{Case $\rDiamFour$:}
		\[
		\infer[\rDiamFour]{\Gamma,\Diamond\varphi \seq \Diamond\psi}
		{\varphi \seq \Diamond\psi}.
		\]
		By the induction hypothesis, every formula in the premise derivation is a $\mathcal{T}$-instance of a subformula of the premise end-sequent $\{\varphi, \Diamond\psi\}$. Because $\varphi$ is an immediate subformula of $\Diamond\varphi$, and $\Diamond\psi$ is trivially a subformula of itself, every formula in the upper derivation is established as a $\mathcal{T}$-instance of a subformula of the conclusion. 
		
		\medskip
		The induction is complete, establishing that every formula in $\pi$ is a $\mathcal{T}$-instance of a subformula of some formula in $\Gamma \cup \{C\}$, which satisfies property (i). 
		
		Finally, we establish property (ii). Substituting terms for variables does not alter the modal depth of a formula. Furthermore, if a formula $F$ is a subformula of $G$, it is necessary that $\mathrm{md}(F) \le \mathrm{md}(G)$. Because property (i) guarantees that every formula in $\pi$ is a $\mathcal{T}$-instance of a subformula of some formula in $\Gamma \cup \{C\}$, its modal depth cannot exceed the maximum modal depth of the formulas in $\Gamma \cup \{C\}$. This establishes property (ii) and completes the proof.
	\end{proof}
	
	\subsection{Craig Interpolation\\}
	\label{subsec:interpolation}
	For a formula $\varphi$, the \emph{predicate vocabulary} $\vocab(\varphi)$ is the set of relation symbols occurring in $\varphi$, and $\consts(\varphi)$ is the set of individual constants occurring in $\varphi$. For a finite context, these notations denote the corresponding unions.
	
	A \emph{retained language} is a triple
	\[
	\Lambda=(\mathcal P,\mathcal C,X),
	\]
	where $\mathcal P$ is a set of relation symbols, $\mathcal C$ a set of constants, and $X$ a set of variables. A formula $B$ is a \emph{$\Lambda$-formula} if
	\[
	\vocab(B)\subseteq\mathcal P,
	\qquad
	\consts(B)\subseteq\mathcal C,
	\qquad
	\fv(B)\subseteq X.
	\]
	
	We define the projection $A\upharpoonright_\Lambda$ while keeping track of variables bound above the current subformula. For a finite set $Y$ of such variables, define $A\upharpoonright_{\Lambda,Y}$ recursively by
	\begin{align*}
		\top\upharpoonright_{\Lambda,Y}&=\top,\\
		(B\wedge C)\upharpoonright_{\Lambda,Y}
		&=(B\upharpoonright_{\Lambda,Y})\wedge(C\upharpoonright_{\Lambda,Y}),\\
		(\Diamond B)\upharpoonright_{\Lambda,Y}
		&=\Diamond(B\upharpoonright_{\Lambda,Y}),\\
		(\forall xB)\upharpoonright_{\Lambda,Y}
		&=\forall x(B\upharpoonright_{\Lambda,Y\cup\{x\}}).
	\end{align*}
	For an atom $S(t_1,\ldots,t_n)$, put
	\[
	S(t_1,\ldots,t_n)\upharpoonright_{\Lambda,Y}
	=
	\begin{cases}
		S(t_1,\ldots,t_n),&\begin{array}{l}
			S\in\mathcal P,\text{ every constant among the }t_i\text{ belongs to }\mathcal C,\\[-1mm]
			\text{and every variable among the }t_i\text{ belongs to }X\cup Y,
		\end{array}\\[3mm]
		\top,&\text{otherwise.}
	\end{cases}
	\]
	Finally, $A\upharpoonright_\Lambda:=A\upharpoonright_{\Lambda,\varnothing}$. The construction is invariant up to $\alpha$-equivalence: renaming a binder also renames its occurrences and the corresponding member of $Y$, so the atomic retention test is unchanged. Moreover, $A\upharpoonright_\Lambda$ is a $\Lambda$-formula. Thus bound variables are retained, while atoms involving a private predicate, a private constant, or a private free parameter are replaced by $\top$.
	
	\begin{lemma}[Contextual Projection]
		\label{lem:projection}
		Let $A\in\mathcal L^+$ and let $\Lambda=(\mathcal P,\mathcal C,X)$ be a retained language. Then:
		\begin{enumerate}
			\item[(i)] $\vdash A\seq A\upharpoonright_\Lambda$;
			\item[(ii)] if every formula of $\Delta$ and the formula $C$ are $\Lambda$-formulas, then
			\[
			\vdash A,\Delta\seq C
			\quad\Longrightarrow\quad
			\vdash A\upharpoonright_\Lambda,\Delta\seq C.
			\]
		\end{enumerate}
	\end{lemma}
	\begin{proof}
		\textbf{(i)} More generally, for every finite set $Y$ of variables we prove
		$\vdash A\seq A\upharpoonright_{\Lambda,Y}$ by structural induction on $A$. If $A=\top$, use $\rTopR$. If an atom is retained, use General Identity; if it is erased, use $\rTopR$. The conjunction and modal cases follow from the induction hypothesis by $\rAndL,\rAndR$ and $\rDiamK$, respectively.
		
		Let $A=\forall xB$. By induction,
		\[
		\vdash B\seq B\upharpoonright_{\Lambda,Y\cup\{x\}}.
		\]
		Choose $y$ fresh from all formulas under consideration and apply Lemma~\ref{lem:subst} with $[y/x]$ to obtain
		\[
		\vdash B[y/x]\seq
		(B\upharpoonright_{\Lambda,Y\cup\{x\}})[y/x].
		\]
		Weakening adds $\forall xB$ to the antecedent; $\rAllL$ then removes the displayed instance $B[y/x]$, and $\rAllR$ gives
		\[
		\vdash \forall xB\seq
		\forall x(B\upharpoonright_{\Lambda,Y\cup\{x\}})
		=(\forall xB)\upharpoonright_{\Lambda,Y}.
		\]
		Taking $Y=\varnothing$ proves (i).
		
		\textbf{(ii)} By soundness and finite completeness, it suffices to prove the semantic entailment
		\[
		A\upharpoonright_\Lambda,\Delta\models C.
		\]
		Let $\mathcal M=(W,R,D,I)$, let $w\in W$, and let $g$ be a $w$-assignment such that
		\[
		\mathcal M,w\models_g A\upharpoonright_\Lambda
		\quad\text{and}\quad
		\mathcal M,w\models_g\Delta.
		\]
		We may first $\alpha$-rename the bound variables of $A$ so that none is a free variable of $A$; this changes neither $A$ nor its projection up to $\alpha$-equivalence.
		
		Let
		\[
		\mathcal C_A=\consts(A)\setminus\mathcal C,
		\qquad
		X_A=\fv(A)\setminus X.
		\]
		Choose pairwise distinct fresh individuals
		\[
		E=\{e_c:c\in\mathcal C_A\}\cup\{e_x:x\in X_A\}.
		\]
		Satisfaction is invariant under isomorphisms of $\QRC$-models,
		by a routine induction on formulas. Hence, because the ambient universe
		$\mathcal U$ is countably infinite and $E$ is finite, transport the
		domains, interpretations, and assignment along an isomorphism to a copy
		disjoint from $E$, and relabel that copy as $\mathcal M$ and its
		transported assignment as $g$. Define an expanded model
		\[
		\mathcal M^+=(W,R,D^+,I^+)
		\]
		by
		\[
		D^+(v)=D(v)\cup E
		\qquad(v\in W).
		\]
		The domains remain nonempty and expanding. Retained constants keep their old interpretations; for $c\in\mathcal C_A$ set $I^+(c)=e_c$; all other constants may keep their old interpretations. For an $n$-ary relation symbol $S$, put
		\[
		I_v^+(S)=
		\begin{cases}
			D^+(v)^n,&S\notin\mathcal P,\\[1mm]
			I_v(S)\cup
			\{\bar d\in D^+(v)^n:\bar d\text{ contains an element of }E\},
			&S\in\mathcal P.
		\end{cases}
		\]
		Finally define $g^+$ by $g^+(x)=e_x$ for $x\in X_A$ and $g^+(x)=g(x)$ otherwise.
		
		We use two generalized structural inductions. First, for every formula $B$ with $\vocab(B)\subseteq\mathcal P$ and $\consts(B)\subseteq\mathcal C$, every world $v$, and every old $v$-assignment $h$,
		\begin{equation}
			\label{eq:fresh-preservation}
			\mathcal M,v\models_hB\quad\Longleftrightarrow\quad\mathcal M^+,v\models_hB.
		\end{equation}
		For right-to-left, every predicate occurring in an atomic subformula of $B$ belongs to $\mathcal P$, and $I_v^+(S)$ agrees with $I_v(S)$ on old tuples; modal witnesses preserve old assignments, and the universal quantifier in $\mathcal M$ ranges over the old domain $D(v)$. For left-to-right, fix $d\in D(v)$. For every $u$ such that $u=v$ or $vRu$, define
		\[
		r_u:D^+(u)\to D(u)
		\]
		by fixing every element of $D(u)$ and sending every element of $E$ to $d$. This is well defined because $D(v)\subseteq D(u)$, and the maps are compatible along accessibility: if $uRu'$, then $r_{u'}$ and $r_u$ agree on $D^+(u)$. Induction simultaneously over these worlds $u$ and expanded assignments proves
		\[
		\mathcal M,u\models_{r_u\circ h^+}B
		\quad\Longrightarrow\quad
		\mathcal M^+,u\models_{h^+}B.
		\]
		The modal case uses compatibility of the retractions at the witness, and the universal case uses
		\[
		(r_u\circ h^+)[r_u(a)/x]=r_u\circ(h^+[a/x]).
		\]
		Taking $u=v$ and $h^+=h$, with $h$ old, gives $r_v\circ h=h$ and establishes the left-to-right implication in \eqref{eq:fresh-preservation}.
		
		Applying \eqref{eq:fresh-preservation} to $g$ and then Coincidence in $\mathcal M^+$ yields
		$\mathcal M^+,w\models_{g^+}A\upharpoonright_\Lambda,\Delta$.
		Second, for every occurrence $B$ of a subformula of $A$, let $Y_B$ be the variables bound above it. Simultaneously for all $u$ with $u=w$ or $wRu$ and all $u$-assignments $h^+$ agreeing with $g^+$ outside $Y_B$, induction on $B$ gives
		\[
		\mathcal M^+,u\models_{h^+}B
		\quad\Longleftrightarrow\quad
		\mathcal M^+,u\models_{h^+}B\upharpoonright_{\Lambda,Y_B}.
		\]
		A retained atom is unchanged. An erased atom is true because its predicate is private, a private constant contributes an element of $E$, or it contains $z\notin X\cup Y_B$; then $z$ is free in $A$ and $h^+(z)=g^+(z)=e_z\in E$. In the modal case, transitivity keeps every witness in the cone of $w$; in the universal case, $h^+[a/x]$ agrees with $g^+$ outside $Y_B\cup\{x\}$. At $B=A$ this yields $\mathcal M^+,w\models_{g^+}A$. Soundness gives $\mathcal M^+,w\models_{g^+}C$; Coincidence and the right-to-left direction of \eqref{eq:fresh-preservation} give $\mathcal M,w\models_gC$. Thus $A\upharpoonright_\Lambda,\Delta\models C$, and Theorem~\ref{thm:finite-completeness} yields the desired derivation.
	\end{proof}
	
	\begin{theorem}[Craig Interpolation]
		\label{thm:interpolation}
		If $\vdash \Gamma_1,\Gamma_2 \seq C$ in $\GQRC$, then there exists a formula $\theta\in\mathcal{L}^+$ such that
		\begin{enumerate}
			\item[(1)] $\vdash \Gamma_1 \seq \theta$;
			\item[(2)] $\vdash \Gamma_2, \theta \seq C$;
			\item[(3)] $\vocab(\theta) \subseteq \vocab(\Gamma_1)\cap \vocab(\Gamma_2\cup\{C\})$;
			\item[(4)] $\consts(\theta) \subseteq \consts(\Gamma_1)\cap \consts(\Gamma_2\cup\{C\})$;
			\item[(5)] $\fv(\theta) \subseteq \fv(\Gamma_1)\cap \fv(\Gamma_2\cup\{C\})$.
		\end{enumerate}
	\end{theorem}
	\begin{proof}
		By Cut Admissibility (Theorem~\ref{thm:cut}), we may use cut-free derivations whenever applying the structural lemmas below.
		Let
		\[
		A=\bigwedge\Gamma_1
		\qquad(\bigwedge\varnothing:=\top)
		\]
		and retain precisely the language of the right-hand side:
		\[
		\Lambda=
		\bigl(
		\vocab(\Gamma_2\cup\{C\}),
		\consts(\Gamma_2\cup\{C\}),
		\fv(\Gamma_2\cup\{C\})
		\bigr).
		\]
		Set
		\[
		\theta=A\upharpoonright_\Lambda.
		\]
		By Lemma~\ref{lem:projection}(i), $A\seq\theta$ is derivable. If $\Gamma_1\neq\varnothing$, repeated applications of height-preserving $\wedge L$-invertibility (Lemma~\ref{lem:andl-inv}) give $\Gamma_1\seq\theta$; if $\Gamma_1=\varnothing$, then $A=\theta=\top$ and (1) is $\rTopR$.
		
		If $\Gamma_1\neq\varnothing$, repeated applications of $\rAndL$ to the assumed derivation of $\Gamma_1,\Gamma_2\seq C$ give
		\[
		A,\Gamma_2\seq C.
		\]
		If $\Gamma_1=\varnothing$, the same sequent is obtained by weakening with $\top$.
		Every formula in $\Gamma_2$ and the formula $C$ is a $\Lambda$-formula. Lemma~\ref{lem:projection}(ii) therefore yields
		\[
		\theta,\Gamma_2\seq C,
		\]
		which is (2).
		
		Finally, $\theta$ is obtained from $A$, so it contains only symbols and free variables occurring in $\Gamma_1$; by construction it retains only those also occurring in $\Gamma_2\cup\{C\}$. This proves (3)--(5).
	\end{proof}
	
	\subsection{Uniform Interpolation\\}
	The preceding projection works independently of the particular retained-language context and conclusion. It therefore yields the following contextual form of left uniform interpolation.
	\begin{theorem}[Contextual Left Uniform Interpolation]
		\label{thm:uniform}
		Let $A\in\mathcal L^+$ and let $\Lambda=(\mathcal P,\mathcal C,X)$ be a retained language. Put
		\[
		A_\Lambda=A\upharpoonright_\Lambda.
		\]
		Then:
		\begin{enumerate}
			\item[(1)] $\vdash A\seq A_\Lambda$;
			\item[(2)] for every finite context $\Delta$ consisting of $\Lambda$-formulas and every $\Lambda$-formula $C$,
			\[
			\vdash A,\Delta\seq C
			\quad\Longrightarrow\quad
			\vdash A_\Lambda,\Delta\seq C.
			\]
		\end{enumerate}
		Equivalently, for every such $\Delta$ and $C$,
		\[
		\vdash A,\Delta\seq C
		\quad\Longleftrightarrow\quad
		\vdash A_\Lambda,\Delta\seq C.
		\]
		Indeed, the reverse implication follows from~(1) by Cut. In particular, with $\Delta=\varnothing$, $A_\Lambda$ is the left uniform interpolant of $A$ in the retained language $\Lambda$.
	\end{theorem}
	\begin{proof}
		Both clauses are exactly Lemma~\ref{lem:projection}. The formula $A_\Lambda$ depends only on $A$ and $\Lambda$, and not on the choice of $\Delta$ or $C$, which is the required uniformity.
	\end{proof}
	
	\begin{remark}
		Theorem~\ref{thm:uniform} strengthens the predicate-only formulation of uniform interpolation: it uniformly eliminates private relation symbols, individual constants, and free parameters. Theorem~\ref{thm:interpolation} is obtained by applying this contextual result to the left part of a split antecedent and retaining the full language of the right part.
		
		The corresponding property of \emph{right} uniform interpolation fails in $\QRC$. Under the standard dual formulation, a right uniform interpolant for a formula $C$ in a retained predicate vocabulary $\mathcal P$ would be a formula $C^{\mathcal P}$ such that
		\[
		\vocab(C^{\mathcal P})\subseteq\mathcal P,
		\qquad
		C^{\mathcal P}\vdash C,
		\]
		and, moreover, for every formula $B$ with $\vocab(B)\subseteq\mathcal P$,
		\[
		B\vdash C
		\quad\Longrightarrow\quad
		B\vdash C^{\mathcal P}.
		\]
		We show that, in general, no formula can satisfy even the entailment requirement $C^{\mathcal P}\vdash C$.
		
		Assume the signature contains distinct unary relation symbols $P$ and $Q$, let
		\[
		C=P(x)\wedge Q(x),
		\qquad
		\mathcal P=\{P\},
		\]
		and suppose that some $C^{\mathcal P}$ satisfies
		\[
		\vocab(C^{\mathcal P})\subseteq\{P\}
		\qquad\text{and}\qquad
		C^{\mathcal P}\vdash P(x)\wedge Q(x).
		\]
		Consider the single-world model $\mathcal M=(\{w\},R,D,I)$ given by
		\[
		R=\{(w,w)\},
		\qquad
		D(w)=\{d\}.
		\]
		Every individual constant denotes $d$, while
		\[
		I_w(P)=\{d\},
		\qquad
		I_w(Q)=\varnothing;
		\]
		all other relation symbols may be interpreted as empty. This is a $\QRC$-model.
		
		For every assignment $g$ at $w$ and every formula $B$ with
		$\vocab(B)\subseteq\{P\}$, one has
		\[
		\mathcal M,w\models_g B.
		\]
		This follows by structural induction. The cases of $\top$ and conjunction are immediate. Every atom is of the form $P(t)$, and every term denotes $d$, so every atom is true. In the universal case the induction hypothesis applies under the unique updated assignment $g[d/x]$, and in the modal case the induction hypothesis at $w$, together with $wRw$, supplies the required witness.
		
		In particular, $\mathcal M,w\models_g C^{\mathcal P}$, whereas
		$\mathcal M,w\not\models_gP(x)\wedge Q(x)$ because
		$I_w(Q)=\varnothing$. Hence
		\[
		C^{\mathcal P}\not\models P(x)\wedge Q(x),
		\]
		and soundness gives
		\[
		C^{\mathcal P}\nvdash P(x)\wedge Q(x),
		\]
		a contradiction. Thus right uniform interpolation fails in general.
	\end{remark}
	
	\section{Computational Complexity}
	\label{sec:complexity}
	In this section, we establish that the decision problem for $\QRC$ is $\mathbf{NP}$-complete over the fixed signatures specified below. We use the standard finite-string encoding of formulas over a fixed countable signature with a polynomial-time decidable presentation, and choose fresh eigenvariables with polynomial-length names. The upper bound applies to every such signature. For the lower bound, we assume that the signature contains a binary relation symbol $R$ and an effectively indexed sequence of pairwise distinct constants $c_0,c_1,\ldots$ such that $c_i$ can be produced in time polynomial in $\log(i+2)$ and has length polynomial in $\log(i+2)$. We bridge the gap between decidability provided by the finite model property and the exact complexity class by leveraging the unique structural constraints of the $\GQRC$ calculus.
	The absence of left-branching rules, such as those for disjunction or implication on the left, is one of the structural reasons why derivable sequents in $\GQRC$ admit polynomial-size cut-free derivations. Furthermore, its modal rules act as strict structural bottlenecks. The context-erasing nature of $\rDiamK$ and $\rDiamFour$ forces the derivation to permanently discard unselected formulas, keeping the proof tree narrow. Conversely, as we will show in the lower bound, this exact bottlenecking mechanism perfectly encodes the non-deterministic ``guessing and verifying'' required for $\mathbf{NP}$-hardness.
	
	\subsection{Upper Bound: NP\\}
	We first prove that the derivability decision problem belongs to $\mathbf{NP}$ by showing that every derivable input sequent has a polynomial-size, polynomial-time checkable cut-free derivation.
	At first glance, quantified modal logics risk exponential blowup due to the interaction between modal contexts and nested universal quantifiers (e.g., expanding a formula like $\forall x \forall y P(x,y)$ with every term in a domain). However, the absence of left-branching rules and the strict context-erasure of its modalities allow us to establish the existence of polynomial-size cut-free derivations despite the potentially large space of redundant instantiations. Consequently, a non-deterministic algorithm can guess a valid derivation of polynomial size and deterministically verify it in polynomial time.
	
	\begin{definition}[Antecedent-occurrence ancestry]
		\label{def:occurrence-ancestry}
		Orient a derivation from its end-sequent toward its axioms. Antecedent occurrences are connected to occurrences in premises as follows. A parametric occurrence is connected to its corresponding occurrence in each premise; at $\rAndR$ it is connected to both copied occurrences. At $\rAndL$, the principal occurrence $A\wedge B$ is connected to the two generated occurrences $A$ and $B$. At $\rAllL$, the principal occurrence $\forall xA$ is connected both to its retained copy $\forall xA$ and to the generated instance $A[t/x]$; the former edge is called a \emph{continuation edge}, the latter a \emph{generation edge}. At a modal inference, the principal occurrence $\Diamond A$ is connected to the singleton occurrence $A$ in the premise, while all parametric context occurrences have no successor. A \emph{descendant} is an occurrence reachable by zero or more of these edges. A \emph{demand} is either one fixed antecedent occurrence of the succedent atom at an $\rId$ leaf (chosen arbitrarily if there are several) or the principal antecedent occurrence of an $\rDiamK$ or $\rDiamFour$ inference.
	\end{definition}
	
	\begin{lemma}[Pruning]
		\label{lem:pruning}
		Let $\pi$ be a normalized cut-free derivation. Suppose an $\rAndL$ inference has no demand descending from either generated conjunct, or an $\rAllL$ inference has no demand descending from its generated instance. Then the end-sequent of $\pi$ has a normalized cut-free derivation with strictly fewer inferences.
	\end{lemma}
	\begin{proof}
		Fix the indicated inference $\iota$ and consider the subderivation above its premise. Mark precisely the descendants, in the sense of Definition~\ref{def:occurrence-ancestry}, of the occurrences generated by $\iota$. In the $\rAndL$ case we also use the undecomposed formula $A\wedge B$ as a \emph{carrier}; in the $\rAllL$ case no extra carrier is needed because the retained occurrence $\forall xA$ already continues into the premise.
		
		We prove by induction on the height of the subderivation above $\iota$ that it can be reconstructed after deleting every marked occurrence, and, in the $\rAndL$ case, carrying one occurrence of $A\wedge B$ parametrically until a modal inference erases it. If the last inference acts only on unmarked occurrences, reconstruct its premise or premises by the induction hypothesis and reapply the same rule. At $\rAndR$, do this in both premises, carrying the same undecomposed formula into both copied contexts. If a left rule acts on a marked occurrence, every occurrence generated or retained from that principal occurrence is marked; after deletion, the reconstructed premise and conclusion coincide, so that left-rule application is omitted. A modal inference cannot select a marked occurrence, because its principal antecedent occurrence would then be a demand, contrary to the hypothesis. If it selects an unmarked occurrence, all marked context occurrences are erased in the original inference; reapply the same modal rule, placing the carrier, when present, in the regenerated conclusion context. At an $\rId$ leaf the distinguished atom is unmarked, so deleting the marked context and adding the carrier still leaves an $\rId$ instance. An $\rTopR$ leaf remains an axiom.
		
		Only the eigenvariable condition for $\rAllR$ requires verification. In the $\rAllL$ pruning case we only delete antecedent occurrences, so no new free variable is introduced. In the $\rAndL$ case, before the carrier is erased, the following invariant holds along every branch:
		\[
		\fv(A\wedge B)\subseteq
		\bigcup\{\fv(D):D\text{ is a currently marked occurrence}\}.
		\]
		It holds initially with equality. It is preserved by parametric propagation and copying; $\rAndL$ preserves the union of free variables of its principal formula; and $\rAllL$ retains its principal universal formula among the marked descendants, so it cannot remove any free variable needed for the inclusion. If an unmarked modal formula is selected, both the marked occurrences and the carrier are erased. Therefore an eigenvariable absent from the original antecedent of an $\rAllR$ inference is also absent after replacement by the carrier. Every reconstructed $\rAllR$ instance is legitimate.
		
		At the bottom of the reconstruction, the premise of $\iota$ has become its conclusion: in the $\rAndL$ case the generated conjuncts have been replaced by the original $A\wedge B$, and in the $\rAllL$ case the generated instance has simply been deleted while the retained universal remains. Hence $\iota$ itself is omitted. No rule, term, eigenvariable, or formula outside the original normalized derivation is introduced, and at least the inference $\iota$ is deleted. No $\rAllR$ inference is removed, so $E$ is unchanged; surviving $\rAllL$ witnesses are unchanged, the carrier already occurred in the original derivation, and modal depth does not increase. Thus all three normalization clauses are preserved, and the resulting cut-free derivation has the same end-sequent and strictly fewer inferences.
	\end{proof}
	
	\begin{lemma}[Polynomial Derivation Size]
		\label{lem:poly-cert}
		If $\vdash \Gamma \seq C$ in $\GQRC$, there exists a cut-free derivation $\pi$ where the total number of inferences is bounded by $\mathcal{O}(N^3)$, with $N = \Vert\Gamma \seq C\Vert$ defined as the total syntactic length of the input sequent.
	\end{lemma}
	\begin{proof}
		By the Subformula Property (Theorem~\ref{thm:subfml}), there exists a normalized cut-free derivation whose instantiation terms and free parameters belong to a pool $\mathcal{T}=T_0^*\cup E$ (bound-variable names are immaterial up to $\alpha$-conversion). Among all such normalized derivations, choose one with the least number of inferences. This class is nonempty by Theorem~\ref{thm:subfml}, and Lemma~\ref{lem:pruning} preserves normalization because its construction introduces no formula, term, or eigenvariable outside the original normalized derivation and leaves every remaining $\rAllL$ witness unchanged. By definition, $|T_0^*|=\max\{1,|T_0|\}\le N$. The set of fresh eigenvariables $E$ is generated exclusively by $\rAllR$ applications, whose number is bounded by $N$ in step~(i) below. Therefore $|\mathcal T|\le 2N$. 
		
		To guarantee the derivation cannot grow exponentially, we globally bound its structural components and the density of left-rule inferences required to complete it:
		
		\begin{enumerate}
			\item[i] Because the strictly positive language $\mathcal{L}^+$ lacks implication ($\to$) and disjunction ($\vee$), $\GQRC$ possesses no left-branching rules. The derivation splits into separate branches \emph{only} at the two-premise right-conjunction rule $\rAndR$. Tracing root-first, left rules leave the succedent unchanged; $\rDiamK$ passes to the immediate scope of its principal diamond; $\rDiamFour$ leaves that diamond occurrence unchanged; and the right rules pass to immediate subformula occurrences. Associate each succedent occurrence with its originating position in the syntax tree of $C$, independently of term instantiation. Left rules and $\rDiamFour$ preserve this position, whereas $\rAllR$ and $\rDiamK$ descend to the unique child position. At $\rAndR$, the two premises enter disjoint child subtrees. Once a branch leaves a position, it cannot return to it. In particular, $\rAndR$ consumes one occurrence of $\wedge$ and sends its two premises to the two disjoint child occurrences, while $\rAllR$ consumes one occurrence of $\forall$. 
			
			Therefore, the total number of $\rAndR$ applications across the entire derivation is bounded by the initial number of $\wedge$ symbols in $C$, and the total number of $\rAllR$ applications is bounded by the initial number of $\forall$ symbols in $C$. Since both quantities are bounded by the total input length $N$, there are globally at most $N$ inferences of $\rAndR$ and $N$ inferences of $\rAllR$. Consequently, the derivation possesses at most $N$ two-premise inferences, which yields a total of at most $N+1$ terminating axioms. The derivation is thus composed of at most $N+1$ continuous branches (where a branch is defined as a sequence of inferences from the end-sequent to an axiom).
			
			\item[ii] For an antecedent $\Delta$, put
			$m(\Delta)=\max\{\mathrm{md}(D):D\in\Delta\}$, with
			$m(\varnothing)=0$. Every non-modal rule leaves $m$ unchanged or
			decreases it. Along any single branch, every modal rule discards
			the unselected context and replaces its principal antecedent
			$\Diamond B$ by $B$, so it strictly decreases $m$. Since the initial modal depth of the end-sequent is bounded by $N$, any single branch can support a sequence of at most $N$ modal inferences. Multiplying the maximum number of branches ($N+1$) by the maximum number of modal rules per branch ($N$), the total number of modal inferences in the entire derivation is bounded above by $N(N+1)$.
			
			\item[iii] By Lemma~\ref{lem:pruning} and the minimality of the chosen derivation, every $\rAndL$ inference has a demand below at least one generated conjunct, and every $\rAllL$ inference has a demand below its generated instance. Assign each left inference to one such demand.
			
			For a fixed demand, follow its unique root-to-demand branch and trace its occurrence ancestry back to an occurrence in the end-sequent; although $\rAndR$ creates two context copies, the fixed demand lies in only one premise, and both copies have the same parent occurrence. Give each occurrence on this lineage the rank equal to the logical complexity of its formula. A $\rAllL$ inference is assigned to this demand only when the lineage passes through its generated instance $B[t/x]$, never merely through the retained occurrence $\forall xB$. Moving root-first from the end-sequent occurrence toward the demand, an assigned $\rAndL$ step passes from $B\wedge D$ to a proper conjunct, an assigned $\rAllL$ step passes from $\forall xB$ to an instance of $B$, and a modal step passes from $\Diamond B$ to $B$. Each such transition lowers rank by at least one, while parametric propagation and $\rAndR$-copying leave it unchanged. Thus modal boundaries cannot reset the measure: the whole lineage is controlled by the complexity of its originating end-sequent occurrence. Since that complexity is at most $N$, at most $N$ left inferences can be assigned to any one demand.
			
			By (i) and (ii), there are at most $N+1$ axioms and $N(N+1)$ modal inferences. The number of demands is therefore at most $(N+1)+N(N+1)=(N+1)^2$, and the total number of left-rule applications is at most $N(N+1)^2$.
		\end{enumerate}
		
		Since every inference belongs to exactly one of these categories, the total number of inferences is bounded above by the sum of their respective bounds. Combining the right-rules (bounded by $2N$) with the left rules (bounded by $N(N+1)^2$), and with the modal rules (bounded by $N^2+N$), and the terminating axioms (bounded by $N+1$), the total number of inferences in the minimal cut-free derivation is bounded above by $N(N+1)^2 + (N^2+N) + 2N + (N+1) = N^3 + 3N^2 + 5N + 1$. This polynomial is $\mathcal{O}(N^3)$. By Lemma~\ref{lem:rename}, a permutation of variables fixing every variable occurring in the end-sequent can send the at most $N+1$ auxiliary variables in $\mathcal T$ that do not occur in the end-sequent to distinct fresh variables with polynomial-length names, preserving normalization and the number of inferences. Under the usual binary indexing of variables, these names can be chosen to have length $O(\log(N+2))$. Because every instantiation draws exclusively from the polynomial pool $\mathcal{T}$ and the language contains no function symbols, every term has depth one and every formula occurrence has polynomial length (linear under the usual token encoding, and at worst quadratic under an explicit bit encoding of symbol names). Along a branch, each left inference adds at most two formula occurrences, so every antecedent contains only polynomially many polynomial-length formulas. Bound variables may be canonically renamed (or represented by de Bruijn indices), so their names also have polynomial bit length. With $O(N^3)$ sequents, the total written size of the derivation is therefore polynomial (a crude bound is $O(N^8)$ under the bit encoding).
	\end{proof}
	\begin{theorem}[NP Membership]
		\label{thm:np-upper}
		For every fixed countable signature with a polynomial-time decidable presentation under the finite-string encoding described above, the derivability problem for $\QRC$ is in $\mathbf{NP}$.
	\end{theorem}
	\begin{proof}
		Given an input sequent $\Gamma \seq C$ of length $N$, a non-deterministic Turing machine can guess a derivation of size polynomial in $N$, as guaranteed by Lemma~\ref{lem:poly-cert}. For each necessary $\rAllL$ application within the derivation, the machine guesses the exact term from the restricted polynomial pool $\mathcal{T}$ required to expose the target subformula. Verification of the derivation against the purely local inference rules of $\GQRC$, including eigenvariable conditions and capture-avoiding substitutions, is then performed deterministically in polynomial time.
	\end{proof}
	\subsection{Lower Bound: NP-Hardness\\}
	To establish the lower bound, we encode the known $\mathbf{NP}$-complete \textsc{Complete Partition} problem into $\GQRC$. The problem asks whether a graph admits a complete, not necessarily proper, vertex colouring with $k$ colours, equivalently, whether it can be mapped onto $K_k$ so that every target edge is the image of a source edge. The maximum such $k$ is the \emph{pseudo-achromatic number}; see Gupta \cite{Gupta1969} and the hardness reference \cite{halldorsson2007complete}. We may restrict to $k\ge3$: instances with $k<3$ are decidable directly, and a polynomial reduction sends each such instance to a fixed yes-instance ($H=K_3,k=3$) or a fixed no-instance (a three-vertex edgeless graph with $k=3$), while leaving instances with $k\ge3$ unchanged.
	Let $V_H$ and $V_K$ be finite vertex sets, and let $E_H$ and $E_K$ be sets of $2$-element subsets of the corresponding vertex sets. Given $H=(V_H,E_H)$ and $K_k=(V_K,E_K)$, the question is whether there is a map $f:V_H\to V_K$ such that, for every $\{a,b\}\in E_K$, some $\{u,v\}\in E_H$ satisfies $\{f(u),f(v)\}=\{a,b\}$. After a polynomial-time relabelling, take
	\[
	V_H=\{1,\ldots,n\},
	\qquad
	V_K=\{0,\ldots,k-1\}.
	\]
	For $k\ge3$ the edge-cover condition already implies that $f$ is surjective. Hence every yes-instance satisfies $k\le n$; if $k>n$, the reduction outputs the fixed no-instance above. We may therefore assume throughout the construction that $3\le k\le n$, even when $k$ is supplied in binary.
	\begin{theorem}[NP-Hardness]
		\label{thm:np-hard}
		For every fixed countable signature satisfying the presentation assumptions above and containing the binary relation symbol $R$ together with the efficiently indexed constants $c_0,c_1,\ldots$, the derivability problem for $\QRC$ is $\mathbf{NP}$-hard.
	\end{theorem}
	\begin{proof}
		To represent undirected edges, encode each one by the conjunction of its two orientations. We use the convention $\bigwedge\varnothing=\top$. 
		
		For the target graph $K_k$, we use the distinct constants $\{c_v \mid v \in V_K\}$ to represent its vertices. For the query graph $H$, we assign a unique variable $x_u$ to each vertex $u \in V_H$. We construct the antecedent $\Gamma_H$ to encode the query graph, universally quantifying over its vertices and wrapping the symmetric edges in a modality. The succedent $C_K$ encodes the target graph:
		\begin{align*}
			\Gamma_H &= \forall x_1 \dots \forall x_n \Diamond \left( \bigwedge_{\{u,v\} \in E_H} \bigl( R(x_u, x_v) \land R(x_v, x_u) \bigr) \right) \\
			C_K &= \Diamond \left( \bigwedge_{\{a,b\} \in E_K} \bigl( R(c_a, c_b) \land R(c_b, c_a) \bigr) \right)
		\end{align*}
		Under the stated encoding, the construction has size
		\[
		O\!\left((|V_H|+|E_H|+k^2)\log^d(|V_H|+k+2)\right)
		\]
		for some fixed integer $d\geq 1$, and is therefore computable in polynomial time. We prove that $\Gamma_H \seq C_K$ is derivable in $\GQRC$ if and only if there exists a valid complete partition $f: V_H \to V_K$.
		
		\medskip\noindent
		$(\Leftarrow)$ \textbf{Constructing the derivation from a valid mapping.} \\
		Assume a valid surjective edge-cover mapping $f$ exists. We construct a cut-free derivation of $\Gamma_H \seq C_K$ root-first.
		
		We begin by applying the left-universal rule $\rAllL$ upwards exactly $|V_H|$ times to the antecedent $\Gamma_H$. For each universally quantified variable $x_u$, we substitute it with the specific constant $c_{f(u)}$ designated by our valid mapping $f$. Because $\rAllL$ retains its principal universal formula, this process also retains the intermediate partially instantiated formulas. It therefore yields a premise of the form
		\[
		\Gamma_H,\;\Pi_f,\;
		\Diamond \left( \bigwedge_{\{u,v\} \in E_H} \bigl( R(c_{f(u)}, c_{f(v)}) \land R(c_{f(v)}, c_{f(u)}) \bigr) \right) \seq C_K,
		\]
		where $\Pi_f$ is the finite multiset of retained partial instances. With the displayed $\Diamond$-formula exposed in the antecedent, we apply the context-erasing modal rule $\rDiamK$ upwards, designating it as principal. This rule strips the outer modalities from both sides and erases the unselected context $\Gamma_H,\Pi_f$, leaving the purely propositional premise:
		\[
		\bigwedge_{\{u,v\} \in E_H} \bigl( R(c_{f(u)}, c_{f(v)}) \land R(c_{f(v)}, c_{f(u)}) \bigr) \seq \bigwedge_{\{a,b\} \in E_K} \bigl( R(c_a, c_b) \land R(c_b, c_a) \bigr)
		\]
		We must now verify that this propositional sequent is provable. We apply the right-conjunction rule $\rAndR$ upwards repeatedly to the succedent. This splits the proof tree into $2|E_K|$ separate branches, where the succedent of each branch is reduced to a single, directed atomic target edge: either $R(c_a, c_b)$ or $R(c_b, c_a)$.
		
		Consider any such branch, requiring us to prove, without loss of generality, $R(c_a, c_b)$. Because our hypothesis guarantees that $f$ is a valid complete partition, we know that for the undirected target edge $\{a,b\}$, there exists an undirected source edge $\{u,v\} \in E_H$ such that $\{f(u), f(v)\} = \{a,b\}$, and $\Gamma_H$ was built precisely out of such edges. 
		
		Our substitution in the antecedent replaced $x_u$ with $c_{f(u)}$ and $x_v$ with $c_{f(v)}$. Since $\{f(u),f(v)\} = \{a,b\}$, after substitution the antecedent contains both $R(c_a,c_b)$ and $R(c_b,c_a)$: one as $R(c_{f(u)},c_{f(v)})$, the other as $R(c_{f(v)},c_{f(u)})$. In particular, the atom $R(c_a,c_b)$ appears among the antecedent conjuncts. We simply apply the left-conjunction rule $\rAndL$ upwards repeatedly to unpack the antecedent until $R(c_a, c_b)$ is exposed, at which point the branch cleanly closes via the $\rId$ axiom. This holds for all branches, completing the derivation.
		
		\medskip\noindent
		$(\Rightarrow)$ \textbf{Extracting a valid mapping from a derivation.} \\
		Conversely, suppose that $\vdash \Gamma_H \seq C_K$. By Theorem~\ref{thm:subfml}, choose a normalized cut-free derivation satisfying the stated term-pool condition. We trace this derivation root-first.
		
		We first analyze the term pool available for substitutions. By the Subformula Property (Theorem~\ref{thm:subfml}), the active term pool is $T_0^*\cup E$. No application of $\rAllR$ can occur, because every succedent subformula is built only from diamonds, conjunctions, and atoms; hence $E=\varnothing$. The active term pool is therefore restricted to the terms explicitly occurring in the end-sequent $\Gamma_H \seq C_K$, all of which belong to $\{c_a\mid a\in V_K\}\cup\{x_u\mid u\in V_H\}$. 
		
		Tracing root-first, the outermost logical structure of $\Gamma_H$ is a block of universal quantifiers. Before the first modal inference, the only applicable rule is $\rAllL$ on the outer universal block or on one of its partial instances. The diamond selected by the first modal inference has a unique occurrence ancestry containing exactly $|V_H|$ generation edges, and therefore determines one witness tuple $(t_u)_{u\in V_H}$. Because the rule $\rAllL$ retains the original quantified formula in the antecedent, it can be applied multiple times, possibly with different terms. Each complete chain of $|V_H|$ applications adds one fully instantiated diamond while the original $\Gamma_H$ and any partial instances are retained. Consequently, after finitely many $\rAllL$ steps, the antecedent may contain several diamond formulas, each corresponding to a different sequence of term substitutions:
		\[
		\Gamma_H,\;\Pi,\; \Diamond A_1,\; \Diamond A_2,\; \dots,\; \Diamond A_m \seq C_K,
		\]
		where $\Pi$ is the multiset of retained partial instances and each $A_i$ is the propositional conjunction $\bigwedge_{\{u,v\}\in E_H} \bigl( R(t_u, t_v) \land R(t_v, t_u) \bigr)$ obtained by instantiating every $x_u$ with some terms $t_u$ drawn from the available pool.
		
		Crucially, no matter how many such diamonds accumulate, the succedent $C_K$ is still a single diamond formula. To expose its propositional content and eventually reach an axiom, the derivation must apply a modal rule. Both $\rDiamK$ and $\rDiamFour$ require selecting exactly one principal diamond formula from the antecedent, permanently erasing all other formulas present there, including the original $\Gamma_H$, the partial instances in $\Pi$, and any other $\Diamond A_j$. Hence, only one of the instantiations, the one chosen as the principal formula, actually influences the remainder of the proof; the others are discarded and irrelevant.
		
		It remains to determine which modal rule is applied. The rule $\rDiamFour$ would produce a premise of the form
		\[
		A_i \seq \Diamond \bigl( \bigwedge_{\{a,b\}\in E_K} \bigl( R(c_a,c_b) \land R(c_b,c_a) \bigr) \bigr).
		\]
		After decomposing $A_i$ by $\rAndL$, the antecedent becomes purely atomic (or consists of $\top$ when $E_H=\varnothing$), while the succedent retains its leading $\Diamond$. At this point no rule of the calculus can be applied:
		\begin{itemize}
			\item Left rules ($\rAndL$, $\rAllL$) cannot affect the $\Diamond$ on the right;
			\item Right rules ($\rAndR$, $\rAllR$) do not match the $\Diamond$-connective;
			\item Modal rules ($\rDiamK$, $\rDiamFour$) require a $\Diamond$-formula in the conclusion's antecedent, which is no longer present;
			\item The axioms $\rId$ and $\rTopR$ are inapplicable because the succedent is neither atomic nor $\top$.
		\end{itemize}
		Thus the branch deadlocks; it can never close. Equivalently, a straightforward induction on the conjunction structure of $A_i$ shows that no cut-free derivation of this premise exists. Therefore, a successful derivation must use $\rDiamK$ instead.
		
		The premise of $\rDiamK$ strips the modality from both sides, yielding a
		purely conjunctive-atomic sequent:
		\[
		\bigwedge_{\{u,v\}\in E_H}
		\bigl(R(t_u,t_v)\land R(t_v,t_u)\bigr)
		\seq
		\bigwedge_{\{a,b\}\in E_K}
		\bigl(R(c_a,c_b)\land R(c_b,c_a)\bigr).
		\]
		A structural induction on any cut-free derivation of a sequent
		$\Delta\seq D$ built only from atoms, $\top$, and conjunction shows
		that every atomic leaf of $D$ occurs as an atomic leaf of some formula
		in $\Delta$: the axioms are immediate, $\rAndR$ uses the induction
		hypothesis for both conjuncts, and $\rAndL$ preserves the antecedent's
		atomic leaves. Hence, for every target atom $R(c_a,c_b)$, some source
		edge $\{u,v\}\in E_H$ contributes either
		$R(t_u,t_v)=R(c_a,c_b)$ or $R(t_v,t_u)=R(c_a,c_b)$. Thus either
		$t_u=c_a$ and $t_v=c_b$, or $t_v=c_a$ and $t_u=c_b$.
		
		We can now explicitly construct the required vertex mapping $f: V_H \to V_K$. For every vertex $u$ such that $t_u=c_v$ for some $v\in V_K$, set $f(u)=v$; assign all remaining vertices arbitrarily. Under this definition, the condition $(t_u = c_a \text{ and } t_v = c_b)$ or $(t_v = c_a \text{ and } t_u = c_b)$ translates directly to the set equality $\{f(u), f(v)\} = \{a, b\}$. Because every atomic leaf of the succedent occurs among the
		atomic leaves of the antecedent, every undirected target edge
		$\{a,b\}\in E_K$ has a source edge $\{u,v\}\in E_H$ satisfying
		$\{f(u),f(v)\}=\{a,b\}$. Therefore the constructed function $f$ is a
		valid complete partition.
	\end{proof}
	\begin{theorem}[Computational Complexity of $\QRC$]
		\label{thm:complexity-final}
		Over every fixed countable signature satisfying the presentation assumptions of this section and containing the binary relation symbol $R$ together with the efficiently indexed constants $c_0,c_1,\ldots$, the decision problem for the Quantified Reflection Calculus $\QRC$ is $\mathbf{NP}$-complete.
	\end{theorem}
	\begin{proof}
		The upper bound is established by Theorem~\ref{thm:np-upper}, and the lower bound is established by Theorem~\ref{thm:np-hard}. Hence the decision problem is $\mathbf{NP}$-complete.
	\end{proof}
	
	\section{Conclusion}
	\label{sec:conclusion}
	In this paper, we have introduced $\GQRC$, the first Gentzen-style sequent calculus for the Quantified Reflection Calculus $\QRC$. Our calculus provides a fully analytic proof theory for a logic that was previously understood primarily through its semantics, its Hilbert-style axiomatization, and its remarkable escape from Vardanyan's undecidability theorem. 
	The primary technical contribution of $\GQRC$ is its resolution of the difficulties associated with proving cut-elimination in the presence of the transitivity axiom ($\Diamond\Diamond\varphi \vdash \Diamond\varphi$). By designing the modal rules $\rDiamK$ and $\rDiamFour$ as strictly \emph{context-erasing} operations, we prevented the accumulation of parametric context across modal boundaries. This structural innovation successfully forced the principal cases of cut-elimination to strictly decrease the inductive derivation measure, allowing us to establish cut-admissibility via a standard lexicographic induction. 
	The analytic control afforded by this cut-elimination theorem yielded a rich harvest of structural and semantic corollaries. We established the subformula property and provided a direct, self-contained proof of the finite model property, guaranteeing decidability. Furthermore, we strengthened the existing interpolation result to a split-context Craig theorem and proved a contextual form of left uniform interpolation that uniformly eliminates private predicate symbols, individual constants, and free parameters.
	Perhaps the most striking consequence of this proof-theoretic framework is our exact characterization of $\QRC$'s computational complexity. We established that the decision problem for $\QRC$ is $\mathbf{NP}$-complete over the fixed countable signatures specified in Section~\ref{sec:complexity}. Within the broader landscape of quantified modal logics, which natively tend toward undecidability (e.g., $\mathrm{QK}$ or full $\mathrm{QPL}(\PA)$) or non-elementary bounds, this result places $\QRC$ in a comparatively low nondeterministic complexity class for a quantified modal logic. The polynomial-certificate proof uses three structural constraints of $\GQRC$: (i) the absence of left-branching rules, (ii) the lack of function symbols, which keeps the active term pool polynomial, and (iii) the context-erasing modal rules, which bound the number of modal steps along each branch. Conversely, our lower-bound reduction from the \textsc{Complete Partition} problem demonstrated how these exact bottlenecks natively model non-deterministic guessing and verification.
	Several promising avenues for future research emerge naturally from this work. The rigid structural bottlenecks and explicit polynomial bounds on cut-free proof trees (Lemma~\ref{lem:poly-cert}) make $\GQRC$ an excellent candidate for developing a root-first proof search algorithm. The polynomial certificate bound suggests implementations for root-first search and polynomial-time verification of candidate derivations, although deterministic proof search may still require exponential time. Furthermore, a natural open question is whether the context-erasing modal architecture of $\GQRC$ can in any way be generalized to polymodal reflection calculi. Managing the interactions between modalities of different strengths without breaking the cut-elimination measure remains a profound proof-theoretic challenge of great significance. In particular, such a generalization would open the way to corresponding arithmetical soundness and completeness questions for polymodal quantified reflection calculi.

\end{document}